\documentclass{article}
\usepackage{amsfonts, amsmath,amssymb}

\newtheorem{teo}{Theorem}

\newtheorem{lemma}{Lemma}

\newtheorem{definition}{Definition}
\newtheorem{remark}{Remark}

\begin{document}

\title{ Nonlinear Hyperbolic-Elliptic Systems in the Bounded Domain}
\author{N.V. Chemetov \\
{\small DCM-FFCLRP/University of Sao Paulo}\\
{\small Ribeirao Preto, Prazil }\\
{\small E-mail: nvchemetov@gmail.com}}
\date{}
\maketitle
\tableofcontents

\bigskip

\textbf{Abstract. }{In the article we study a hyperbolic-elliptic system of
PDE. The system can describe two different physical phenomena: 1st one is
the motion of magnetic vortices in the II-type superconductor and 2nd one \
is the collective motion of cells. Motivated by real physics, we consider
this system with boundary conditions, describing the flux of vortices (and
cells, respectively) through the boundary of the domain. We prove the global
solvability of this problem.\ To show the solvability result we use a
"viscous" parabolic-elliptic system. Since the viscous solutions do not have
a compactness property, we justify the limit transition on a vanishing
viscosity, using a kinetic formulation of our problem. As the final result
of all considerations we have solved a very important question related with
a so-called "boundary layer problem", showing the strong convergence of the
viscous solutions to the solution of our hyperbolic-elliptic system. }



\section{Formulation of the problem}

\label{sec1}\setcounter{equation}{0}

In this paper we consider the system of nonlinear hyperbolic and elliptic
equations, which describes the evolution of \ two unknowns $\omega =\omega
(t,\mathbf{x})$ and $h=h(t,\mathbf{x}):$ 
\begin{eqnarray}
\omega _{t}+\mathrm{div}(g(\omega )\mathbf{v}) &=&0,\quad \quad \quad 
\mathbf{v}=-\nabla h,  \label{sec1eq1} \\
-\Delta h+h &=&\omega \quad \quad \text{in}\quad \Omega _{T}:=(0,T)\times
\Omega ,  \label{sec1eq3}
\end{eqnarray}%
where $g=g(\omega )$ is a nonlinear function of $\omega $ and $\Omega $\ is
a bounded domain in $\mathbb{R}^{2}.$ The system is closed by the boundary
condition for $h=h(t,\mathbf{x})$ on the boundary $\Gamma $\ of the domain $%
\Omega $ 
\begin{equation}
h=a\quad \quad \quad \text{on}\quad \Gamma _{T}:=(0,T)\times \Gamma 
\label{sec1eq5}
\end{equation}%
and by the initial condition for $\omega $%
\begin{equation}
\omega |_{t=0}=\omega _{0}\quad \quad \quad \text{in}\quad \Omega .
\label{sec1eq4}
\end{equation}%
If the characteristics of the hyperbolic equation \eqref{sec1eq1} on the
boundary $\Gamma $\ are directed into the domain $\Omega ,$ then on this
inflow part of $\Gamma $ we need an additional boundary condition for $%
\omega .$ The characteristics on $\Gamma $ are given by the equation $\frac{d%
\mathbf{X}}{dt}=\mathbf{u}(\mathbf{X}(t),t)$ with $\mathbf{u}=g^{\prime
}(\omega )\mathbf{v,}$ such that for some $t_{0}\in \left[ 0,T\right] :$ $\ 
\mathbf{X}(t_{0})=\mathbf{x}_{0}\in \Gamma .$ Therefore denoting the
external normal at $\mathbf{x\in }\Gamma $ by $\mathbf{n}=\mathbf{n}(\mathbf{%
x}),$\ the inflow part of $\Gamma $ is described as $\Gamma _{T}^{\mathbf{-}%
}:=\left\{ (t,\mathbf{x})\in \Gamma _{T}:\quad g^{\prime }(\omega )\left( 
\mathbf{v}\cdot \mathbf{n}\right) (t,\mathbf{x})<0\right\} ,$ where we give
the following boundary condition%
\begin{equation}
\omega =b(t,\mathbf{x},\frac{\partial h}{\partial \mathbf{n}})\quad \text{on}%
\quad \Gamma _{T}^{\mathbf{-}}.  \label{sec1eq6}
\end{equation}

The above system is a two-dimensional reduction of the three dimensional
mean field model for the motion of superconductive vortices, derived in \cite%
{Chap-2}, \cite{Chap-3}. In the system the unknown $h$ is a magnetic field,
applied to the superconductor sample,\ the unknown $\omega $ is a magnetic
vortex density and the function $g$\ \ is defined as $g(\omega ):=|\omega |.$
During the last 15 years this superconductive model has been extensively
studied and a review of these investigations can be found in the articles 
\cite{ant-che-2007}, \cite{Chemet}. In the last two articles the solvability
of our system has been shown but under the assumption of the positivity for
the vortex density $\omega ,$ implying that $g(\omega )\equiv \omega .$
Hence a nonlinear \ hyperbolic equation \eqref{sec1eq1} transforms to a
linear transport equation, that significantly simplifies the analysis of the
system \eqref{sec1eq1}-\eqref{sec1eq6}. Taking into account of the nonlinear
hyperbolicity of \eqref{sec1eq1} and an explosion behavior of the vortices
(see the experimental results given on the web-site:
www.fys.uio.no/super/results/sv/index.html) \ it is more interest to show
the solvability of our system, admitting the explosion behavior.

The same type system has appeared in a model of Keller-Segel to describe
transport effects of the collective motion of cells (usually bacteria or
amoebae), that are attracted by a chemical substance and are able to emit it
(see the review article \cite{hor} and references therein). The unknown $%
\omega $ represents the cell density, the unknown $h$ is the concentration
of chemoattractant, inducing a drift force\bigskip . The function $g(\omega
) $ takes the form $g(\omega ):=\omega (1-\omega ).$

Let us note that our system does not permit to derive $BV$-boundedness of $%
\omega ,$ hence in our case we can not apply the Kruzkov's approach, which
normally is very powerful in the study of hyperbolic equations. The first
positive result of weak solvability for the Keller-Segel's model in the
unbounded domain $\Omega \equiv \mathbb{R}^{2}$ has been obtained by
Perthame and Dalibard \cite{perthame1}, using a beautiful idea: it is enough
to investigate properties of a \textit{linear transport} equation, obtained
from a kinetic reformulation of the system. The same idea has been suggested
firstly by Plotnikov in \cite{plot}, \cite{plot1} (also for a problem in the
unbounded domain $\Omega \equiv \mathbb{R}^{2}).$

In this article we study the solvability of \eqref{sec1eq1}-\eqref{sec1eq6}
without any simplifications and use the method of \cite{perthame1}-\cite%
{plot1}, developing it to a \emph{bounded} domain $\Omega .$ We consider the
problem \eqref{sec1eq1}-\eqref{sec1eq6} for the two models, considering that
1) $g(\omega ):=|\omega |$ or 2) $g(\omega ):=\omega (1-\omega ).$ We will
not separate the proof of the solvability for these two models (Theorem \ref%
{teo2sec1} and Theorem \ref{teo2sec1 copy(1)}), since the only difference is
in obtaining the $L_{\infty }$-estimate for $\omega $ (see the section \ref%
{sec7.2}). Let us note that the $L_{\infty }-$ boundedness for $\omega $ in
two models presents completely two different problems. This estimate for the
Keller-Segel model follows directly from the maximum principle, but for the
mean field model it is possible to obtain just from energy type estimates.
The last one is a main difficulty of the boundary layer theory.

\textit{We assume \ that }$\Omega $\textit{\ is a bounded domain \ of} $%
\mathbb{R}^{2}$\textit{\ having the boundary }$\Gamma $\textit{\ of \ }$C^{2}
$\textit{-smooth and}%
\begin{equation}
\omega _{0}\in L_{\infty }(\Omega ),\quad a\in L_{\infty }(0,T;W_{\infty
}^{2}(\Gamma )),  \label{reg1}
\end{equation}%
\textit{and the function }$b=b(\cdot ,\cdot ,z)$\textit{\ is a continuous
function on }$z\in \mathbb{R},$\textit{\ such that }%
\begin{equation}
|b(t,\mathbf{x},z)|\leqslant b_{0}(t,\mathbf{x})+b_{1}|z|^{\varkappa }\quad 
\text{\textit{for a.e.}}\quad (t,\mathbf{x})\in \Gamma _{T}  \label{reg2}
\end{equation}%
\textit{with some constants }$\varkappa \in (0,1),$\textit{\ }$%
b_{1}\geqslant 0$\textit{\ and a positive \ function }$b_{0}\in L_{\infty
}(\Gamma _{T}).$

Since \eqref{sec1eq1} \ is a nonlinear hyperbolic equation, the solution of %
\eqref{sec1eq1} has to be considered as entropy one.

\begin{definition}
\label{def1} We say that the pair of functions$\quad \omega \in L_{\infty
}(\Omega _{T})$ and $h\in L_{\infty }(0,T;W_{p}^{2}(\Omega ))$\textit{\ for }%
$\forall p<\infty \quad $is a weak solution of the problem \eqref{sec1eq1}-%
\eqref{sec1eq6} if the pair $\left\{ \omega ,h\right\} $ satisfies:

1) the inequality 
\begin{eqnarray}
&&\int_{\Omega _{T}}\left\{ |\omega -\xi |\varphi _{t}+sign(\omega -\xi )
\left[ (g(\omega )-g(\xi ))\,(\mathbf{v}\cdot \nabla \varphi )+g(\xi
)\,(h-\omega )\,\varphi \right] \right\} \ dtd\mathbf{x}  \notag \\
&&\qquad \qquad \qquad \ +\int_{\Omega }|\omega _{0}-\xi |\varphi (0,\cdot
)\,d\mathbf{x}+M(\mathbf{v})\int_{\Gamma _{T}}|b(\cdot ,\cdot ,\frac{%
\partial h}{\partial \mathbf{n}})-\xi |\,\varphi \ dtd\mathbf{x}\geqslant
0\qquad   \label{fraca}
\end{eqnarray}%
for any $\xi \in \mathbb{R}$ and any non-negative $\varphi \in C^{\infty }(%
\overline{\Omega }_{T})$, such that $\varphi (\cdot ,T)=0.$ Here $M(\mathbf{v%
}):=\emph{K}||\mathbf{v}||_{L_{\infty }(\Omega _{T})}$ with $\ \emph{K:=}%
||g^{\prime }||_{L_{\infty }(\mathbb{R})};$

2) the functions $h$ and $\mathbf{\omega }$ fulfill the equation %
\eqref{sec1eq3} a.e. in $\Omega _{T}$ and the boundary condition %
\eqref{sec1eq5} a.e. on $\Gamma _{T}.$
\end{definition}

We refer to Lemmas 7.24 and 7.34 of the book \cite{malek} for a discussion
of the fulfillment of the boundary condition \eqref{sec1eq6} in the sense of
integral inequality \eqref{fraca}. Let us just note that we can take that$\
K\equiv 1$ \ in both cases $g=|\omega |$ and $g=\omega -\omega ^{2},$
accounting the definition (\ref{const}).

Our first result is the following theorem.

\begin{teo}
\label{teo2sec1} Let $g(\omega )=|\omega |.$ If the data $a,\,\,b,\,\,\omega
_{0}$ satisfy (\ref{reg1})-(\ref{reg2}), then there exists at least one weak
solution $\{\omega ,h\}$ of the problem \eqref{sec1eq1}-\eqref{sec1eq6}.
\end{teo}

\begin{remark}
1) The above theorem can be generalized on the case when \textit{\ }$\Omega $%
\textit{\ }is a bounded domain \ of $\mathbb{R}^{n}$ with $n>2.$ But taking
into account that the system \ \eqref{sec1eq1}-\eqref{sec1eq6} is the
2-dimensional reduction of the 3D mean field model and the 3D model has a
more complex structure then our system \ \eqref{sec1eq1}-\eqref{sec1eq6}, we
do not see a reason to improve the theorem for the n-dimensional case with $%
n>2;$

2) In the articles \cite{Chap-2}-\cite{Chap-3} on the inflow boundary $%
\Gamma _{T}^{\mathbf{-}}$ \ the following boundary conditions have been
suggested 
\begin{equation}
-(\mathbf{v}\cdot \mathbf{n})\text{ }\omega =\,|a|\,\frac{B(|\nabla h|,t,%
\mathbf{x})}{|\nabla h|}\,\frac{\partial h}{\partial \mathbf{n}}\quad \quad 
\text{or}\quad \quad -(\mathbf{v}\cdot \mathbf{n})\text{ }\omega
=|a|\,B(|\nabla h|,t,\mathbf{x})\,\frac{\partial h}{\partial \mathbf{n}}
\label{flux}
\end{equation}%
where $a$ is a given function, depending on physical properties of the
superconductor and $B:=\max (|\nabla h|-J,0)$ with $J=J(t,\mathbf{x})>0$
being a so-called function of the nucleation of the vortices on the
boundary. These conditions are equivalent to \textit{\eqref{sec1eq6} with \ }%
$\varkappa =1$ in \eqref{reg2}. But it is not possible to obtain the
existence result when $\varkappa =1$ in \eqref{reg2}, to understand it we
refer the reader to see the deduction of the estimate \eqref{est-2}. In fact
we propose that $B:=\max (|\nabla h|-J,0)^{\varkappa }$ in the 2nd boundary
condition of \eqref{flux}, keeping the main physical effect: a flux
nucleation of vortices at the inflow boundary $\Gamma _{T}^{\mathbf{-}}.$

3)The existence result depends strongly on the $L_{\infty }-$ boundedness of 
$\omega .$ The boundedness follows from the energy type estimate \eqref{ww},
which permits to obtain also $L_{p}-$ boundedness of $\omega $ with $\
p<\infty $. It is open question: Is the existence result valid for the class
of solutions having only $L_{p}-$ boundedness of $\omega $ with some $%
p<\infty ?$
\end{remark}

\bigskip

Our second result is the following theorem.

\begin{teo}
\label{teo2sec1 copy(1)} Let \ $g(\omega )=\omega -\omega ^{2}.$ Let the
data $a,\,\,\,\,\omega _{0},$ $b$ satisfy (\ref{reg1})-(\ref{reg2}), such
that 
\begin{equation}
\quad \quad \quad \,0\leqslant \omega _{0}\leqslant 1\quad \quad \quad \text{%
\textit{on }}\Omega \quad \quad \text{\textit{and}}\quad \quad 0\leqslant
b_{0}(t,\mathbf{x})\leqslant 1,\quad b_{1}=0\quad \text{\textit{on}}\quad
\Gamma _{T}.  \label{wa1}
\end{equation}%
Then there exists at least one weak solution $\{\omega ,h\}$ of the problem %
\eqref{sec1eq1}-\eqref{sec1eq6}, such that 
\begin{equation}
\,\,0\leqslant \omega \leqslant 1\quad \text{\textit{a.e. in \quad }}\Omega
_{T}.  \label{wa2}
\end{equation}
\end{teo}

\bigskip

\bigskip \textit{The article can be divided on two parts:}

\textit{1}$^{\mathit{st}}$\textit{\ part is } the sections \ref{sec3} and %
\ref{sec32}, where we explain principal ideas of the construction of the
solution for the \ problem \eqref{sec1eq1}-\eqref{sec1eq6}.

In the section \ref{sec3}, we consider a "viscous" parabolic-elliptic
problem of \eqref{sec1eq1}-\eqref{sec1eq6} (see the system \eqref{pe1}-%
\eqref{pe2}) and formulate the solvability result (Lemma \ref{lemma1}) of
the problem \eqref{pe1}-\eqref{pe2}. Since viscous solutions have not a
property of strong convergence on a vanishing viscosity, we write our
viscous problem in a kinetic formulation and pass to a weak limit in
obtained kinetic inequalities \eqref{2.60}-\eqref{2.61} \ and equation %
\eqref{2.8}.

In the section \ref{sec32} we investigate properties of \ this kinetic
formulation. We show that the solution of the kinetic equation has initial
and boundary conditions with values equal only to 0 and 1 (Lemma \ref{lemma2}%
). Since the kinetic equation is a \textit{linear transport }one, the
solution of it has zero and unit values too. This implies that the weak
limit of the viscous solutions is \ the solution of our original problem %
\eqref{sec1eq1}-\eqref{sec1eq6}. As a result of all these considerations we
will obtain a very important result for the theory of boundary layers.

\begin{teo}
\label{theorem3} Let $\{\omega _{\nu },h_{\nu }\}$ be the solutions of \ the
viscous problem \eqref{pe1}-\eqref{pe2} for $\nu >0.$ There exists a
subsequence of $\{\omega _{\nu },h_{\nu }\},$ such that 
\begin{equation*}
\omega _{\nu }\rightarrow \omega \quad \quad \text{\textit{strongly in }}%
L_{\infty }(0,T;L_{p}(\Omega ))\quad \text{\textit{and}}\quad h_{\nu
}\rightharpoonup h\quad \quad \text{strongly in $L_{\infty
}(0,T;W_{p}^{1}(\Omega ))$}
\end{equation*}%
for any $p<\infty ,$ where $\{\omega ,h\}$ is a weak solution of the problem %
\eqref{sec1eq1}-\eqref{sec1eq6}.
\end{teo}

As we know this theorem is the first positive convergence result of viscous
solutions to a nonviscous solution with the boundary conditions, which admit
the flow through the boundary of the domain. We do not know a positive
convergence result even for a more "simple" case: Navier-Stokes and Euler
equations (see a review in \cite{cons}, \cite{Kel}).

\textit{2}$^{\mathit{nd}}$\textit{\ part of the article\ is }the section \ref%
{appendix}, where we prove all technical results, formulated in the sections %
\ref{sec3} and \ref{sec32}.
 
This article is a version of the article ”Nonlinear Hyperbolic-Elliptic Systems in the Bounded Domain ” published in \cite{chem2}.  After the publication of this paper,  similar PDE systems, being the nonlinear hyperbolic-elliptic systems have been published  such as  \cite{chem1}-\cite{chem5}, where the Kruzkov method of entropies have been generalized by the development of a so-called ''kinetic method''. We also mention the articles , where the last kinetic model has been applied to porous media models. 

Since the boundary layer problem has importance in the study of motion of fluids we also refer to  \cite{ACC25} and \cite{C}-\cite{CC6}. The stochastic perturbation for a similar system can be seen also in  \cite{ACC21}, \cite{CC0}.

\section{"Viscous" Problem and Kinetic Formulation}

\label{sec3}\setcounter{equation}{0}

Let us consider $\nu \in (0,1)$ and approximate our data $a,$ $\omega _{0},$ 
$g,$ $b,$ $b_{0}$\ by smooth functions $a^{\nu },$ $\omega _{0}^{\nu },$ $%
g^{\nu },$ $b^{\nu },$ $b_{0}^{\nu }$ \ (the construction of this
approximation is given in the section \ref{appendix}). In this section for
simplicity of notations, we drop the parameter {$\nu $} and indicate the
dependence of functions and constants on {$\nu $}${,}$ where it will be
necessary.

For a fixed $\ \nu $ we study the "viscous" approximation of the problem %
\eqref{sec1eq1}-\eqref{sec1eq6} 
\begin{equation}
\left\{ 
\begin{array}{l}
\omega _{t}+\mbox{div}(g(\omega )\mathbf{v})={\ \nu }\Delta \omega ,\quad
\quad \mathbf{v}=-\nabla h\text{ }\quad \text{in }\quad \Omega _{T}, \\ 
\\ 
\nu \frac{\partial \omega }{\partial \mathbf{n}}+M(\mathbf{v}(t))(\omega
-b(\cdot ,\text{\textbf{$\cdot $}},\frac{\partial h}{\partial \mathbf{n}}%
))=0\quad \quad \text{on}\quad \Gamma _{T}, \\ 
\\ 
\omega |_{t=0}=\omega _{0}\quad \quad \text{in}\quad \Omega%
\end{array}%
\right.  \label{pe1}
\end{equation}%
with $M(\mathbf{v}(t)):=K\left\Vert \mathbf{v}(t,\cdot )\right\Vert
_{L_{\infty }(\Omega )}$ and 
\begin{equation}
\left\{ 
\begin{array}{ll}
-\Delta h+h=\omega & \text{in}\quad \Omega _{T}, \\ 
h=a & \text{on}\quad \Gamma _{T}.%
\end{array}%
\right.  \label{pe2}
\end{equation}

\bigskip

The following result is shown in the section \ref{sec7}.

\begin{lemma}
\label{lemma1} There exists at least one solution $\{\omega ,h\}$ of the
problem \eqref{pe1}-\eqref{pe2}, such that%
\begin{equation*}
\omega \in C^{1+\alpha /2,2+\alpha }(\overline{\Omega }_{T}),\quad \quad
h\in C^{\alpha /2,2+\alpha }(\overline{\Omega }_{T})
\end{equation*}%
for some $\alpha \in (0,1).$ \ The pair $\left\{ \omega ,h\right\} $ $\ $%
fulfills the estimates%
\begin{eqnarray}
||\omega ||_{L_{\infty }(\Omega _{T})} &\leqslant &C,  \label{est-2} \\
\left\Vert h\right\Vert _{L_{\infty }(0,T;W_{p}^{2}(\Omega )\cap C^{1+\alpha
}(\overline{\Omega }))} &\leqslant &C,\quad \forall p<\infty ,\quad \forall
\alpha \in \lbrack 0,1)  \label{est3} \\
\sqrt{\nu }||\nabla \omega ||_{L_{2}(\Omega _{T})} &\leqslant &C,
\label{est22} \\
||\partial _{t}\nabla \left( h-h_{a}\right) ||_{L_{2}(\Omega _{T})}
&\leqslant &C  \label{est1}
\end{eqnarray}%
For the trace value $\frac{\partial h}{\partial \mathbf{n}}$ on $\Gamma ,$
we have 
\begin{equation}
\frac{\partial \left( h-h_{a}\right) }{\partial \mathbf{n}}\in \mathcal{P}%
:=L_{2}\big(0,T;W_{2}^{\frac{1}{2}}(\Gamma )\big)\cap W_{2}^{\frac{1}{2}}%
\big(0,T;L_{2}(\Gamma )\big),  \notag
\end{equation}%
such that 
\begin{equation}
||\frac{\partial \left( h-h_{a}\right) }{\partial \mathbf{n}}||_{\mathcal{P}%
}\!\leqslant C.  \label{est2}
\end{equation}%
Here \bigskip $h_{a}$\ is the solution of the system%
\begin{equation}
\left\{ 
\begin{array}{ll}
-\Delta h_{a}+h_{a}=0 & \text{in}\quad \Omega _{T}, \\ 
h_{a}=a & \text{on}\quad \Gamma _{T}.%
\end{array}%
\right.  \label{haha}
\end{equation}%
The constants $C$ are independent of $\nu .$
\end{lemma}

\bigskip

\bigskip Now we may take the limit on $\nu \rightarrow 0$ in the system (\ref%
{pe1})-(\ref{pe2}). However, the estimates of Lemma \ref{lemma1} do not
guarantee the strong convergence of a subsequence of $\left\{ \omega _{\nu
}\right\} $ and do not permit to pass to the limit in the nonlinear term of (%
\ref{pe1}). By this reason the main objective in the following
considerations is to show the strong convergence of $\left\{ \omega _{\nu
}\right\} .$ \ To do it we use the kinetic formulation of (\ref{pe1}).

Let $(\eta ,q)$ be an entropy pair, i.e. $\eta =\eta (s)\in C^{2}(\mathbb{R}%
) $ is a convex function and $q^{\prime }(s):=\eta ^{\prime }(s)g^{\prime
}(s). $ Let $(\eta ,q)$ satisfy the conditions 
\begin{equation}
\eta \geqslant 0\quad \text{and}\quad |q|\leqslant K\eta \quad \text{on}%
\quad \mathbb{R}.  \label{ineqeta}
\end{equation}%
(also take into account (\ref{const})). From the equation of (\ref{pe1}), we
have 
\begin{equation*}
\partial _{t}\eta +\mathrm{div}(\mathbf{v}q)+(\omega -h)(g\eta ^{\prime
}-q)=\nu \eta ^{\prime }\Delta \omega
\end{equation*}%
with $\eta ,$ $\eta ^{\prime },$ $q,$ $g$ \ being the functions of $\omega .$%
\ Let us multiply this equality by $\varphi \in C^{\infty }(\overline{\Omega 
}_{T}),$ such that $\varphi |_{t=T}=0$ and integrate it over $\Omega _{T}.$
Using that $\eta (b)=\eta (\omega )+\eta ^{\prime }(\omega )(b-\omega )+%
\frac{\eta ^{\prime \prime }(r)}{2}(b-\omega )^{2}$ for some $r=r(\mathbf{x}%
,t)$ with values between $\omega (\mathbf{x},t)$ and $b(\mathbf{x},t,\frac{%
\partial h}{\partial \mathbf{n}}(\mathbf{x},t)),$ we get that $\omega $
fulfills the entropy type equality 
\begin{eqnarray}
\int_{\Omega _{T}}\{\eta \varphi _{t} &+&q\,(\mathbf{v}\cdot \nabla \varphi
)+(g\eta ^{\prime }-q)\,(h-\omega )\,\,\varphi \}dtd\mathbf{x}+\int_{\Omega
}\eta (\omega _{0})\varphi (0,\text{\textbf{$\cdot $}})\,\,d\mathbf{x} 
\notag \\
&+&\int_{\Gamma _{T}}M(\mathbf{v}(t))\eta (b(\cdot ,\cdot ,\frac{\partial h}{%
\partial \mathbf{n}}))\varphi \ dtd\mathbf{x}-\nu \int_{\Omega _{T}}\eta
^{\prime }\nabla \omega \nabla \varphi \ dtd\mathbf{x}=m_{\nu ,\eta
}(\varphi )\quad \quad  \label{2}
\end{eqnarray}%
with 
\begin{equation*}
m_{\nu ,\eta }(\varphi ):=\int_{\Omega _{T}}\nu \eta ^{\prime \prime
}|\nabla \omega |^{2}\varphi \ dtd\mathbf{x}+\int_{\Gamma _{T}}\{(\mathbf{%
v\cdot n})q+M(\mathbf{v}(t))\eta +\frac{1}{2}M(\mathbf{v}(t))\eta ^{\prime
\prime }(r)(b-\omega )^{2}\}\varphi \ dtd\mathbf{x}.
\end{equation*}%
\medskip Observe that $m_{\nu ,\eta }(\varphi )\geqslant 0$ for $\ \varphi
\geqslant 0$ as a consequence of (\ref{ineqeta}).

Let us denote by $\left\vert s\right\vert _{\pm }:=\max (\pm s,0)$ and
consider $\eta (s):=|s|_{+}.$ Since this $\eta \notin C^{2}(\mathbb{R})$ \
we first take a $smooth$ convex function $\eta _{\varepsilon },$ such that $%
\eta _{\varepsilon }(s)\rightarrow \eta (s)$ strongly in $C_{loc}(\mathbb{R)}
$ for $\varepsilon \rightarrow 0.$ The entropy pair $(\eta _{\varepsilon
}(\omega -\xi ),$ $q_{\varepsilon }(\omega -\xi ))$ with an arbitrary $\xi
\in \mathbb{R}$ satisfies the equality (\ref{2}). Hence using that 
\begin{equation*}
\left( \eta _{\varepsilon }(\omega -\xi ),\,q_{\varepsilon }(\omega -\xi
)\right) \rightarrow \left( |\omega -\xi |_{+},\,\mathrm{sign}_{+}(\omega
-\xi )(g(\omega )-g(\xi ))\right) =:(\eta ,q)
\end{equation*}%
strongly in $C(\overline{\Omega }_{T})$ as $\varepsilon \rightarrow 0$, we
infer that (\ref{2}) is valid for this entropy pair $(\eta ,q).$ Therefore,
the function 
\begin{equation*}
f_{\nu }(\xi ,t,\mathbf{x}):=\mathrm{sign}_{+}(\omega _{\nu }-\xi ):=%
\begin{cases}
1, & \text{if }\omega _{\nu }(t,\mathbf{x})>\xi , \\ 
0, & \text{if }\omega _{\nu }(t,\mathbf{x})\leqslant \xi%
\end{cases}%
\end{equation*}%
satisfies%
\begin{align}
& \int\limits_{\Omega _{T}}\Bigg\{\int\limits_{\xi }^{+\infty }f_{\nu }(s,t,%
\mathbf{x})\{\varphi _{t}+g^{\prime }(s)(\mathbf{v}\cdot \nabla \varphi
)\}\,ds+g(\xi )f_{\nu }(\xi ,t,\mathbf{x})(h_{\nu }-\omega _{\nu })\varphi %
\Bigg\}\,dtd\mathbf{x}  \notag \\
& +\int\limits_{\Omega }|\omega _{0,\nu }-\xi |_{+}\varphi (0,\mathbf{x})\,d%
\mathbf{x}\ +\int\limits_{\Gamma _{T}}M(\mathbf{v}_{\nu }(t))|b_{\nu }(t,%
\mathbf{x},\frac{\partial h_{\nu }}{\partial n})-\xi |\,_{+}\varphi \ dtd%
\mathbf{x}\,  \notag \\
& -\nu \int\limits_{\Omega _{T}}f_{\nu }(\xi ,t,\mathbf{x})\nabla \omega
_{\nu }\nabla \varphi \,dtd\mathbf{x}=m_{\nu ,|\omega _{\nu }-\xi
|_{+}}(\varphi )  \label{in1}
\end{align}%
\ with 
\begin{eqnarray}
m_{\nu ,|\omega -\xi |_{+}}(\varphi ):= &&\int_{\Omega _{T}}\nu \delta (\xi
=\omega _{\nu }(t,\mathbf{x}))\ |\nabla \omega _{\nu }|^{2}\varphi \ dtd%
\mathbf{x}  \notag \\
&+&\int_{\Gamma _{T}}\{(\mathbf{v}_{\nu }\cdot \mathbf{n})\ f_{\nu }(\xi ,t,%
\mathbf{x})(g(\omega _{\nu })-g(\xi ))+M(\mathbf{v}_{\nu }(t))\ |\omega
_{\nu }-\xi |_{+}\qquad  \notag \\
&+&\frac{1}{2}M(\mathbf{v}_{\nu }(t))\ \delta (\xi =s_{\nu }(t,\mathbf{x}%
))(b_{\nu }-\omega _{\nu })^{2}\varphi \}\ dtd\mathbf{x,}\qquad  \label{me}
\end{eqnarray}%
such that $m_{\nu ,|\omega -\xi |_{+}}(\varphi )\geqslant 0$ for $\varphi
\geqslant 0.$ Here $\delta (s)$ is the Dirac function.

Choosing $\varphi (t,\mathbf{x}):=1-1{_{\varepsilon }}(t-(T-\varepsilon ))$
in \eqref{in1} and taking the limit on $\varepsilon \rightarrow 0,$ with the
help of (\ref{est-2}), (\ref{est3}) we get 
\begin{equation}
m_{\nu ,|\omega -\xi |_{+}}(1)\mathbf{\leqslant }\int\limits_{\Omega
}|\omega _{0,\nu }-\xi |_{+}\ d\mathbf{x}+\int\limits_{\Gamma _{T}}M(\mathbf{%
v}_{\nu }(t))\ |b_{\nu }(t,\mathbf{x},\frac{\partial h_{\nu }}{\partial n}%
)-\xi |\,_{+}\,\ dtd\mathbf{x}+Cg(\xi )\mathbf{.}  \label{mes}
\end{equation}%
In view of Riesz representation theorem (see \cite{evans}) the measure $%
m_{\nu }(\xi ,t,\mathbf{x}):=m_{\nu ,|\omega (t,\mathbf{x})-\xi |_{+}}$ is
well defined on $\mathbb{R\times }\overline{\Omega }_{T}.$ Let us remark
properties of $m_{\nu }:$ \ 

1) $m_{\nu }\in BV(\mathbb{R},w-\mathcal{M}^{+}(\overline{\Omega }_{T}))$ in
view of the left part of the equality (\ref{in1})$.$ Therefore there exist
left and right limits of $m_{\nu }(\xi ,\cdot ,\cdot )$ at any $\xi \in 
\mathbb{R}.$ Here \ $\mathcal{M}^{+}(\Omega _{T})$ \ is the space of \
bounded \ non-negative Radon measures on $\Omega _{T}$ and$\ \ w-\mathcal{M}%
^{+}(\Omega _{T})$ denotes \ $\mathcal{M}^{+}(\Omega _{T})$ equipped with
the weak topology.

2) $m_{\nu }(\xi ,\cdot ,\cdot )=0$ on $\Omega _{T}$ and for $|\xi |>R:=%
\underset{\nu }{\sup }||\omega _{\nu }||_{L_{\infty }(\Omega )},$ since $%
m_{\nu }(\xi ,\cdot ,\cdot )=\nu \delta (\xi =\omega _{\nu }(\cdot ,\cdot
))\ |\nabla \omega _{\nu }|^{2}$ in any subdomain of $\Omega _{T}.$

\bigskip

If we take $\eta :=|\omega -\xi |_{-}$ \ for $\xi \in \mathbb{R}$\ in %
\eqref{2}, by a similar way as \ (\ref{in1}) has been deduced, we show 
\begin{align}
& \int\limits_{\Omega _{T}}\Bigg\{\int\limits_{-\infty }^{\xi }(1-f_{\nu
}(s,t,\mathbf{x}))\{\varphi _{t}+g^{\prime }(s)(\mathbf{v}\cdot \nabla
\varphi )\}\,ds+g(\xi )(f_{\nu }(\xi ,t,\mathbf{x})-1)(h_{\nu }-\omega _{\nu
})\varphi \Bigg\}\,dtd\mathbf{x}  \notag \\
& +\int\limits_{\Omega }|\omega _{0,\nu }-\xi |_{-}\varphi (0,\cdot )\,d%
\mathbf{x}\,+\int\limits_{\Gamma _{T}}M(\mathbf{v}_{\nu }(t))\ |b_{\nu }(t,%
\mathbf{x},\frac{\partial h_{\nu }}{\partial n})-\xi |\,_{-}\,\varphi \ dtd%
\mathbf{x}\qquad  \notag \\
& -\nu \int\limits_{\Omega _{T}}(f_{\nu }(\xi ,t,\mathbf{x})-1)\nabla \omega
_{\nu }\nabla \varphi \,dtd\mathbf{x}\geqslant 0\quad \quad \text{for}\quad 
\text{ }\varphi \geqslant 0.  \label{in2}
\end{align}%
Now if we take $\varphi :=\frac{\partial \psi }{\partial \xi }$ with $\psi
(\xi ,t,\mathbf{x})\in C_{0}^{\infty }(\mathbb{R}\times \Omega _{T})$ in %
\eqref{in1}, then after integration by parts on the parameter $\xi ,$ \ we
derive that $f_{\nu }$ satisfies the following identity%
\begin{align}
& \int\limits_{\mathbb{R}\times \Omega _{T}}f_{\nu }\left\{ \psi
_{t}+g^{\prime }(\xi )(\mathbf{v}\cdot \nabla _{\mathbf{x}}\psi )+g(\xi
)(h_{\nu }-\omega _{\nu })\psi _{\xi }^{\prime }\right\}  \notag \\
& \qquad \ \qquad \ \ -\nu \nabla \omega _{\nu }\nabla \psi _{\xi }^{\prime
}\}\,d\xi dtd\mathbf{x}=\int\limits_{\mathbb{R}\times \Omega _{T}}\psi _{\xi
}^{\prime }\,m_{\nu }(\xi ,t,\mathbf{x})\ d\xi dtd\mathbf{x}.  \label{eq1}
\end{align}%
Let us collect properties of $f_{\nu },$ $\omega _{\nu },$ $m_{\nu },$ used
in the sequel 
\begin{equation}
\begin{cases}
0\leqslant f_{\nu }\leqslant 1\quad \quad \text{ }\quad \text{ a.e. in }%
\mathbb{R}\times \Omega _{T}\quad \quad \text{and}\quad \text{ }\frac{%
\partial f_{\nu }}{\partial \xi }\leqslant 0\quad \text{ }\quad \text{ }%
\quad \text{ in \ }\mathcal{D}^{\prime }(\mathbb{R})\text{, } \\ 
f_{\nu }(\xi ,\mathbf{\cdot },\mathbf{\cdot })=%
\begin{cases}
0\quad \text{ for}\quad \text{ }\xi >R, \\ 
1\quad \text{ for}\quad \text{ }\xi <-R,%
\end{cases}%
\quad \quad \text{and}\quad m_{\nu }(\xi ,\mathbf{\cdot },\mathbf{\cdot }%
)=0\quad \text{for}\quad |\xi |>R, \\ 
\omega _{\nu }(t,\mathbf{x})=\int\limits_{0}^{+\infty }f_{\nu }(s,t,\mathbf{x%
})\,ds-\int\limits_{-\infty }^{0}(1-f_{\nu }(s,t,\mathbf{x}))\,ds, \\ 
m_{\nu }\in BV(\mathbb{R},w-\mathcal{M}^{+}(\overline{\Omega }_{T}))\quad 
\text{satisfying \ (\ref{mes}).}%
\end{cases}
\label{eq2.5}
\end{equation}

\section{Limit Transition on $\protect\nu \rightarrow 0$}

\label{sec32}\setcounter{equation}{0}

In this section we prove the strong convergence of a subsequence $\left\{
\omega _{\nu }\right\} $ to $\omega ,$ that implies the solvability of our
main problem \eqref{sec1eq1}-\eqref{sec1eq6}.

In view of Lemma \ref{lemma1} and \eqref{eq2.5}, there exists a subsequence
of $\{\omega _{\nu },f_{\nu },h_{\nu },m_{\nu }\},$ such that 
\begin{align*}
h_{\nu }& \rightharpoonup h\quad \text{weakly -- $\ast $ in $L_{\infty
}(0,T;W_{p}^{2}(\Omega ))\quad \forall p<\infty ,$}\quad \quad \\
\frac{\partial h_{\nu }}{\partial \mathbf{n}}& \rightarrow \frac{\partial h}{%
\partial \mathbf{n}}\quad \text{strongly in }L_{2}(\Gamma _{T}), \\
\mathbf{v}_{\nu }& \rightarrow \mathbf{v}\quad \text{strongly in $L_{\infty
}(0,T;L^{2}(\Omega ))$}, \\
\omega _{\nu }& \rightharpoonup \omega \quad \text{weakly -- $\ast $ in $%
L_{\infty }(\Omega _{T})$},\quad \quad \quad \quad \\
f_{\nu },\,\,\omega _{\nu }f_{\nu }& \rightharpoonup \,f,\,\,\rho \quad 
\text{weakly -- $\ast $ in $L_{\infty }(\mathbb{R}\times \Omega _{T}),$} \\
\nu \nabla \omega _{\nu }& \rightarrow 0\quad \text{strongly in $%
L_{2}(\Omega _{T})$},\quad \quad \quad \quad \\
m_{\nu }& \rightharpoonup m\quad \text{weakly \ in }\mathcal{M}^{+}(\mathbb{%
R\times }\overline{\Omega }_{T}).
\end{align*}%
Moreover $\omega ,f,m$ \ fulfill the following properties 
\begin{equation}
\begin{cases}
0\leqslant f\leqslant 1\quad \quad \text{ a.e. in }\mathbb{R}\times \Omega
_{T}\quad \quad \text{ and}\quad \quad \text{ }\frac{\partial f}{\partial
\xi }\leqslant 0\quad \quad \text{ in \ }\mathcal{D}^{\prime }(\mathbb{R})%
\text{, }\quad \text{ }\quad \\ 
f(\xi ,\mathbf{\cdot },\mathbf{\cdot })=%
\begin{cases}
0\quad \text{ for}\quad \text{ }\xi >R, \\ 
1\quad \text{ for}\quad \text{ }\xi <-R%
\end{cases}%
\quad \text{and}\quad m(\xi ,\mathbf{\cdot },\mathbf{\cdot })=0\quad \text{%
for}\quad |\xi |>R, \\ 
\omega (t,\mathbf{x})=\int\limits_{0}^{+\infty }f(s,t,\mathbf{x}%
)\,ds-\int\limits_{-\infty }^{0}(1-f(s,t,\mathbf{x}))\,ds, \\ 
m\in BV(\mathbb{R},w-\mathcal{M}^{+}(\overline{\Omega }_{T})),\quad \quad 
\text{such that }\forall \xi \in \mathbb{R}:\text{ \ \ \ }\int\limits_{%
\overline{\Omega }_{T}}m(\xi ,t,\mathbf{x})\,dtd\mathbf{x}<C.%
\end{cases}
\label{eq2.88}
\end{equation}%
Obviously, the functions $\omega ,f,h,\mathbf{v}$ satisfy \eqref{sec1eq3}-%
\eqref{sec1eq5}, the estimates \eqref{est-2}, \eqref{est3} and the integral
relations 
\begin{align}
& \int\limits_{\Omega _{T}}\int\limits_{\xi }^{+\infty }f(s,t,\mathbf{x}%
)\left\{ \varphi _{t}+g^{\prime }(s)(\mathbf{v}\cdot \nabla \varphi
)\right\} \,ds+g(\xi )(fh-\rho )\varphi \,dt\,d\mathbf{x}  \notag \\
& \qquad \ +\int\limits_{\Omega }|\omega _{0}-\xi |_{+}\varphi (0,\cdot )\,d%
\mathbf{x}+\int\limits_{\Gamma _{T}}M(\mathbf{v}(t))|b(t,\mathbf{x},\frac{%
\partial h}{\partial n})-\xi |_{+}\varphi \,dt\,d\mathbf{x}\geqslant 0,
\label{2.60}
\end{align}%
\begin{align}
& \int\limits_{\Omega _{T}}\int\limits_{-\infty }^{\xi }(1-f(s,t,\mathbf{x}%
))\left\{ \varphi _{t}+g^{\prime }(s)(\mathbf{v}\cdot \nabla \varphi
)\right\} \,ds+g(\xi )\left\{ (f-1)(h-\omega )+(f\omega -\rho )\right\}
\varphi \,dt\,d\mathbf{x}  \notag \\
& \qquad +\int\limits_{\Omega }|\omega _{0}-\xi |_{-}\varphi (0,\cdot )\,d%
\mathbf{x}+\int\limits_{\Gamma _{T}}M(\mathbf{v}(t))|b(t,\mathbf{x},\frac{%
\partial h}{\partial n})-\xi |_{-}\varphi \,dt\,d\mathbf{x}\geqslant 0
\label{2.61}
\end{align}%
with $M(\mathbf{v}(t))=\emph{K}\left\Vert \mathbf{v}(t,\cdot )\right\Vert
_{L_{\infty }(\Omega _{T})}$ and%
\begin{equation}
\int\limits_{\mathbb{R}\times \Omega _{T}}f\{\partial _{t}\psi +g^{\prime
}(\xi )(\mathbf{v}\cdot \nabla _{\mathbf{x}}\psi )\}+g(\xi )(fh-\rho )\psi
_{\xi }^{\prime }\,\,d\xi \,dt\,d\mathbf{x}=\int\limits_{\mathbb{R}\times
\Omega _{T}}\psi _{\xi }^{\prime }\,m\,d\xi dtd\mathbf{x}.  \label{2.8}
\end{equation}

\vspace{1pt}

The equality \eqref{2.8} guarantees the existence of traces for $f$ \ at the
time $t=0$ and on the boundary $\Gamma .$ Due to the inequalities %
\eqref{2.60}, \eqref{2.61}, the traces have crucial properties formulated in
the following lemma.

\begin{lemma}
\label{lemma2}

\begin{enumerate}
\item[1)] There exists the trace $f^{0}=f^{0}(\xi ,\mathbf{x})$ at the time
moment $t=0$ \ \ for the function $f$, such that 
\begin{equation}
\begin{cases}
f^{0}=\lim\limits_{\varepsilon \rightarrow 0^{+}}\frac{1}{\varepsilon }%
\int_{0}^{\varepsilon }f(\cdot ,t,\mathbf{\cdot })\,dt\text{\ \ \ \ \ a.e. \
on \ \ }\mathbb{R}\times \Omega , \\ 
0\leqslant f^{0}\leqslant 1\text{\ }\quad \text{a.e. \ in \ \ }\mathbb{R}%
\times \Omega \text{\ \ \ \ \ \ and\ \ \ \ \ \ }(f^{0})_{\xi }^{\prime
}\leqslant 0\quad \text{in \ }\mathcal{D}^{\prime }(\mathbb{R}\times \Omega
).%
\end{cases}
\label{2.10}
\end{equation}%
There exist non-negative functions $m_{+}^{0}=m_{+}^{0}(\xi ,\mathbf{x}),$ $%
m_{-}^{0}=m_{-}^{0}(\xi ,\mathbf{x})\in W_{\infty }^{1}(\mathbb{R},L_{\infty
}(\Omega )),$ such that 
\begin{equation}
\begin{cases}
f^{0}=\mathrm{sign}_{+}(\omega _{0}-\xi )+\partial _{\xi }m_{+}^{0},\text{\
\ \ \ \ }\quad \text{\ \ \ \ \ \ \ }\quad \lim\limits_{\xi \rightarrow
+\infty }m_{+}^{0}=0\text{\ }\quad \text{a.e. in }\Omega ; \\ 
f^{0}-1=\mathrm{sign}_{+}(\omega _{0}-\xi )-1+\partial _{\xi }m_{-}^{0},%
\text{\ \ \ \ \ }\quad \lim\limits_{\xi \rightarrow -\infty }m_{-}^{0}=0%
\text{\ }\quad \text{a.e. in }\Omega .%
\end{cases}
\label{2.13}
\end{equation}

\item[2)] There exists the trace $f^{\Gamma }=f^{\Gamma }(\xi ,t,\mathbf{x})$
on the boundary $\Gamma _{T}$ \ for the function $f,$ such that 
\begin{equation}
\begin{cases}
f^{\Gamma }(\cdot ,\cdot ,\mathbf{x})=\lim\limits_{\varepsilon \rightarrow
0^{+}}\frac{1}{\varepsilon }\int_{0}^{\varepsilon }f(\cdot ,\cdot ,\mathbf{x}%
-s\,\mathbf{n}(\mathbf{x}))\,ds\text{ \ on }\mathbb{R}\times \Gamma _{T}:\ \
g^{\prime }(\xi )(\mathbf{v}\cdot \mathbf{n})\neq 0; \\ 
0\leqslant f^{\Gamma }\leqslant 1\text{\ }\quad \text{a.e. \ on \ \ }\mathbb{%
R}\times \Gamma _{T}\text{\ }\quad \text{\ and\ }\quad (f^{\Gamma })_{\xi
}^{\prime }\leqslant 0\text{\ }\quad \text{in \ }\mathcal{D}^{\prime }(%
\mathbb{R}\times \Gamma _{T}).%
\end{cases}
\label{eq2.15}
\end{equation}%
There exist non-negative functions $m_{+}^{\Gamma }(\xi ,t,\mathbf{x})$, $%
m_{-}^{\Gamma }(\xi ,t,\mathbf{x})\in W_{\infty }^{1}(\mathbb{R},L_{\infty
}(\Gamma _{T})),$ such that 
\begin{equation}
\begin{cases}
g^{\prime }(\xi )(\mathbf{v}\cdot \mathbf{n})f^{\Gamma }=-M\mathrm{sign}%
_{+}(b(\cdot ,\mathbf{\cdot },\frac{\partial h}{\partial \mathbf{n}})-\xi
)-\partial _{\xi }m_{+}^{\Gamma }, \\ 
\text{\ \ \ \ \ }\lim\limits_{\xi \rightarrow +\infty }m_{+}^{\Gamma }=0%
\text{\ }\quad \text{a.e. on }\Gamma _{T}; \\ 
g^{\prime }(\xi )(\mathbf{v}\cdot \mathbf{n})(f^{\Gamma }-1)=M\left[ 1-%
\mathrm{sign}_{+}(b(\cdot ,\mathbf{\cdot },\frac{\partial h}{\partial 
\mathbf{n}})-\xi )\right] -\partial _{\xi }m_{-}^{\Gamma }, \\ 
\text{\ \ \ \ \ }\lim\limits_{\xi \rightarrow -\infty }m_{-}^{\Gamma }=0%
\text{\ }\quad \text{a.e. on }\Gamma _{T}.\text{\ \ \ \ \ }\quad%
\end{cases}
\label{2.16}
\end{equation}
\end{enumerate}
\end{lemma}

\vspace{1pt}\qquad

Now let us show two properties for the function $\rho .$ Observe that $%
(\omega _{\nu }-\xi )f_{\nu }(\xi )=\int_{\xi }^{+\infty }f_{\nu }(s)ds.$
Hence taking $\nu \rightarrow 0,$ we derive that $\rho (\xi )=\xi f(\xi
)+\int_{\xi }^{\infty }f(s)\,ds$ \ \ \ a.e. in $\Omega _{T},$ which is 
\begin{equation}
\rho _{\xi }^{\prime }=\xi f_{\xi }^{\prime }\text{ }\quad \text{ }\quad 
\text{ }\quad \text{ }\quad \text{in \ }\mathcal{D}^{\prime }(\mathbb{R}%
\times \Omega _{T})  \label{eq2.20}
\end{equation}
too. Since $\rho =f\omega $ for $|\xi |>R,$ then from \ \eqref{eq2.88} we
have%
\begin{align*}
\rho (\xi )-\omega f(\xi )& =\xi \,f(\xi )+\int_{\xi }^{+\infty
}f(s)\,ds-\omega f(\xi ) \\
& =f(\xi )\left[ \int_{-R}^{\xi }(1-f(s))ds\right] +(1-f(\xi ))\left[
\int_{\xi }^{R}f(s)ds\right]
\end{align*}%
and accounting that $f$ is decreasing on $\xi ,$ we get 
\begin{equation}
|\rho -\omega f|\leqslant 2R\ f\ (1-f)\quad \text{ }\quad \text{ }\quad 
\text{ a.e. in }\mathbb{R}\times \Omega _{T},  \label{eq2.21}
\end{equation}%
wher $R$ is defined in \eqref{eq2.5}.

The equation \eqref{2.8} \ is written as 
\begin{equation}
\partial _{t}f+\mathrm{div}_{\mathbf{x}}(g^{\prime }(\xi )\,\mathbf{v}%
\,f)+(g(\xi )(fh-\rho ))_{\xi }^{\prime }=\partial _{\xi }m\quad \text{ in \ 
}\mathcal{D}^{\prime }(\mathbb{R\times }\Omega _{T}).  \label{eq2.9}
\end{equation}%
Because of \eqref{eq2.20}, $g\in W_{\infty }^{1}(\mathbb{R}),$ $\ g^{\prime
}\in BV(\mathbb{R})$ \ and $\mathbf{v},h\in L_{\infty }(0,T,W_{p}^{1}(\Omega
)),$ we can apply the renormalization theorem to the left part of %
\eqref{eq2.9} \ \ (see for instance Theorem 4.3 in \cite{lellis}) and get,
that the function $F:=f(1-f)$ satisfies%
\begin{align}
& F_{t}+\mathrm{div}_{\mathbf{x}}(g^{\prime }(\xi )\mathbf{v}F)+(g(\xi
)(h-\xi )F)_{\xi }^{\prime }+  \notag \\
& +(g(\xi )(\xi -\omega ))_{\xi }^{\prime }F+g^{\prime }(\xi )(1-2f)\left[
f\omega -\rho \right] \leqslant 0\quad \text{ in \ }\mathcal{D}^{\prime }(%
\mathbb{R\times }\Omega _{T}).  \label{eqq1}
\end{align}%
It means that \eqref{eq2.9} is regularized on a parameter $\theta ,$
multiplied by $(1-2f^{\theta })$ ($f^{\theta }$ being the regularization of $%
f)$\ and finally taken the limit on $\theta \rightarrow 0.$ \ The inequality
in \eqref{eqq1} follows from the relation $\int\limits_{\mathbb{R}}\frac{%
\partial m^{\theta }}{\partial \xi }(1-2f^{\theta })\,d\xi =2\int\limits_{%
\mathbb{R}}m^{\theta }\frac{\,\partial f^{\theta }}{\partial \xi }\,d\xi
\leqslant 0$\ (see the properties \eqref{eq2.88}). The derived inequality %
\eqref{eqq1} implies the following result.

\begin{lemma}
\label{lemma33} For a.e. $t\mathbb{\in }(0,T)$ the inequality holds 
\begin{eqnarray}
\int_{\mathbb{R}\times \Omega }F\ d\xi d{\mathbf{x}}(t) &\leqslant &C\{\int_{%
\mathbb{R}\times \Omega }f^{0}-\left( f^{0}\right) ^{2}\ d\xi d{\mathbf{x}} 
\notag \\
&+&\int\limits_{\mathbb{R}\times \Gamma _{T}}{|g^{\prime }(\xi )(\mathbf{v}%
\cdot \mathbf{n})|}_{-}\,\,(f^{\Gamma }-\left( f^{\Gamma }\right)
^{2})\,\,d\xi dtd\mathbf{x\}}.  \label{uha}
\end{eqnarray}
\end{lemma}

Moreover we can show the following result.

\begin{lemma}
\label{lemma3}The traces $f^{0}$ and $f^{\Gamma }$ satisfy 
\begin{equation*}
\int\limits_{\mathbb{R}\times \Omega }f^{0}-(f^{0})^{2}\,d\xi \,d\mathbf{x}=0%
\text{ \ \ \ and \ \ }\int\limits_{\mathbb{R}\times \Gamma _{T}}{|g^{\prime
}(\xi )(\mathbf{v}\cdot \mathbf{n})|}_{-}\,\,(f^{\Gamma }-\left( f^{\Gamma
}\right) ^{2})\,\,d\xi dtd\mathbf{x}=0.
\end{equation*}
\end{lemma}

By Lemmas \ref{lemma33} and \ref{lemma3}, we obtain that $F=0$. Hence $f$
takes only the values $0$ and $1$ almost everywhere. Since $f$ is monotone
decreasing on $\xi ,$ \ there exists some function $u=u(t,\mathbf{x}),$ such
that 
\begin{equation*}
f(\xi ,t,\mathbf{x})=\mathrm{sign}_{+}(u(t,\mathbf{x})-\xi ).
\end{equation*}%
Therefore we have the $\ast -$weakly convergence in $L_{\infty }(\Omega
_{T}) $ 
\begin{equation*}
\ |\omega _{\nu }|_{+}=\int\limits_{0}^{+\infty }f_{\nu }(s,\mathbf{\cdot },%
\mathbf{\cdot })\ d\xi ,\quad |\omega _{\nu }|_{-}=\int\limits_{-\infty
}^{0}(1-f_{\nu }(s,\mathbf{\cdot },\mathbf{\cdot }))\ d\xi \rightharpoonup
|u|_{+},\quad |u|_{-},
\end{equation*}%
that implies $u=\omega $ (see \eqref{eq2.88}), as a consequence of it we
derive the strong convergence of $\left\{ \omega _{\nu }\right\} $\vspace{1pt%
} to $\omega $ and $\rho =\mathrm{sign}_{+}(\omega -\xi )\ \omega .$ Taking
the sum of the inequalities \eqref{2.60}, \eqref{2.61} \ we see that $\omega 
$ satisfies \eqref{fraca}, that ends the proof of the solvability of our
problem \eqref{sec1eq1}-\eqref{sec1eq6}.

\section{Proofs of Technical Results}

\label{appendix} \setcounter{equation}{0}

\subsection{Proof of Lemma \protect\ref{lemma1}}

\label{sec7}

\bigskip

Let us formulate a well-known result from the approximation theory.

\begin{lemma}
\label{ap1} Let $Q\subseteq \mathbb{R}^{n},\quad n\geqslant 1$ be an open
set and $s\in \lbrack \,1,+\infty ]$. Then for any $\phi \in L_{s}\,(Q)$
there exist functions $\phi ^{\nu }\in C^{\infty }(Q),$ satisfying the
following properties: 
\begin{equation}
||\phi ^{\nu }||_{L_{s}(Q)}\leqslant C||\phi ||_{L_{s}(Q)}\,\quad \,\quad
and\,\quad \quad ||\phi ^{\nu }{-}\phi ||_{L_{r}(Q)}{\longrightarrow 0,}\;%
\mbox{when }\;\nu \rightarrow 0,  \label{ap33}
\end{equation}%
for $r:=s,$ if $s<\infty $ and for any $r<\infty $, if $s=\infty .$ The
constant $C$ is independent of \ the parameter $\nu >0.$
\end{lemma}

\bigskip

Using Lemma \ref{ap1}, we approximate our data $a$, $\omega _{0}$, $b_{0},$ $%
g$ and $b_{z}:=b(\cdot ,\cdot ,z(\cdot ,\cdot ))$ for any $z=z(t,\mathbf{x})$
by $C^{\infty }$ - functions $a^{\nu }$, $\omega _{0}^{\nu }$, $b_{0}^{\nu
}, $ $g^{\nu }$ and $b_{z}^{\nu },$ satisfying the relations \eqref{ap33} in
corresponding spaces $L_{s}(Q)$, defined by the regularity conditions %
\eqref{reg1} - \eqref{reg2}. We can assume that any derivatives of these
smooth approximations are bounded by constants, depending only on $\nu >0.$
For simplicity of notations, we omit the parameter {$\nu $} and indicate the
dependence of functions and constants on {$\nu $}${,}$ where it will be
necessary.

We show the solvability of the problem \eqref{pe1}-\eqref{pe2}, using Leray
- Schauder's fixed point theorem (see \cite{gilbarg}, p. 286, Theorem 11.6).
For it we consider the following problem, depending also on an auxiliary
parameter $\lambda \in \lbrack 0,1]:$

\begin{equation}
\left\{ 
\begin{array}{l}
\omega _{t}+\mbox{div}(\lambda g(\omega )\mathbf{v})={\nu }\Delta \omega
,\;\;\;\mathbf{v}:=-\nabla h\quad \quad \text{\ in}\quad \text{\ }\Omega
_{T}, \\ 
\nu \frac{\partial \omega }{\partial \mathbf{n}}+\lambda M(\mathbf{v}%
(t))(\omega -b(\cdot ,\text{\textbf{$\cdot $}},\frac{\partial h}{\partial 
\mathbf{n}}))=0\quad \quad \text{\ on}\quad \text{\ }\Gamma _{T}, \\ 
\omega |_{t=0}=\lambda \omega _{0}\quad \quad \text{\ in}\quad \text{\ }%
\Omega%
\end{array}%
\right.  \label{h1}
\end{equation}%
with $M(\mathbf{v}(t)):=K\left\Vert \mathbf{v}(t,\cdot )\right\Vert
_{L_{\infty }(\Omega )}$ and

\begin{equation}
\left\{ 
\begin{array}{ll}
-\Delta h+h=\omega & \text{in}\quad \Omega _{T}, \\ 
h=a & \text{on}\quad \text{\ }\Gamma _{T}.%
\end{array}%
\right.  \label{h2}
\end{equation}%
To apply Leray - Schauder's fixed point argument, first we assume the
existence of a solution of \eqref{h1}-\eqref{h2} and deduce the a priori
estimates \eqref{est-2}-\eqref{est3} which do not depend on $\lambda $ and $%
\nu .$ Below until the end of the section \ref{sec7} \ \ all constants $C$ \
will be independent\textit{\ }of $\lambda $ and $\nu .$

\subsubsection{\protect\vspace{1pt} Proof of a priori estimates \eqref{est-2}%
-\eqref{est3}}

\label{sec7.2}

Since the case $\lambda =0$ \ is trivial, in the sequel\ we consider $%
\lambda >0.$ \ Let us note that the solution $\omega $ of \eqref{h1}-%
\eqref{h2} fulfills a similar identity as \eqref{2}. Namely, for any entropy
pair $(\eta ,q)$ satisfying the conditions \eqref{ineqeta}, the following
equality is valid 
\begin{eqnarray}
\int_{\Omega _{T}}\{\eta \varphi _{t} &+&\lambda q\,(\mathbf{v}\cdot \nabla
\varphi )+\lambda (g\eta ^{\prime }-q)\,(h-\omega )\,\,\varphi \}dtd\mathbf{x%
}+\int_{\Omega }\eta (\lambda \omega _{0})\varphi (0,\text{\textbf{$\cdot $}}%
)\,\,d\mathbf{x}\qquad  \notag \\
&+&\lambda \int_{\Gamma _{T}}M(\mathbf{v}(t))\eta (b(\cdot ,\cdot ,\frac{%
\partial h}{\partial \mathbf{n}}))\varphi \ dtd\mathbf{x}-\nu \int_{\Omega
_{T}}\eta ^{\prime }\nabla \omega \nabla \varphi \ dtd\mathbf{x}=m_{\nu
,\eta }(\varphi )\qquad  \label{2222}
\end{eqnarray}%
with 
\begin{eqnarray*}
m_{\nu ,\eta }(\varphi ):= &&\int_{\Omega _{T}}\nu \eta ^{\prime \prime
}|\nabla \omega |^{2}\varphi \ dtd\mathbf{x}+\lambda \int_{\Gamma _{T}}\{(%
\mathbf{v\cdot n})q+M(\mathbf{v}(t))\eta \\
&+&\frac{1}{2}M(\mathbf{v}(t))\eta ^{\prime \prime }(r)(b-\omega
)^{2}\}\varphi \ dtd\mathbf{x}\geqslant \int_{\Omega _{T}}\nu \eta ^{\prime
\prime }|\nabla \omega |^{2}\varphi \ dtd\mathbf{x}
\end{eqnarray*}%
\medskip for any positive $\varphi \in C^{\infty }(\overline{\Omega }_{T}),$
such that $\varphi |_{t=T}=0.$ Here $r=r(\mathbf{x},t)$ is a function having
values between $\omega (\mathbf{x},t)$ and $b(\mathbf{x},t,\frac{\partial h}{%
\partial \mathbf{n}}(\mathbf{x},t)).$ Taking $\varphi (t,\mathbf{x}):=1-1{%
_{\varepsilon }}(t-t_{0})$ for arbitrary $t_{0}\in (0,T)$ and passing to the
limit on $\varepsilon \rightarrow 0$ in \eqref{2222}, we derive%
\begin{align}
& \int_{\Omega }\eta {\ }d\mathbf{x}(t_{0})+\int_{0}^{t_{0}}\int_{\Omega
}\{\nu \eta ^{\prime \prime }|\nabla \omega |^{2}+\lambda \lbrack g\eta
^{\prime }-q]\,(\omega -h)\}{\ }dtd\mathbf{x}  \notag \\
& \qquad \qquad \qquad \mathbf{\leqslant }\int_{\Omega }\eta (\lambda \omega
_{0})\,\,d\mathbf{x}+\int_{0}^{t_{0}}M(\mathbf{v}(t))\int_{\Gamma }\eta
(\lambda b(\cdot ,\cdot ,\frac{\partial h}{\partial \mathbf{n}}))\,\,dtd%
\mathbf{x}  \label{ww}
\end{align}

To derive $L_{\infty }$ - boundedness for the solution $\omega $ of %
\eqref{h1}-\eqref{h2} we consider separately two different cases: when $%
g=|\omega |$ and when $g(\omega )=\omega -\omega ^{2}.$

\textbf{1}$^{\mathbf{st}}$\textbf{\ case:} \ when $g=|\omega |.$ \ Let us
deduce a priori estimate (\ref{est-2}). Let $\eta :=|\omega {}|^{k}$ and $%
q:=sign(\omega {})|\omega {}|^{k}$\ be chosen for any $k\geqslant 2$ in %
\eqref{ww}, that gives the inequality 
\begin{align}
||\omega ||_{L_{k}(\Omega )}^{k}(t_{0})& +\int_{0}^{t_{0}}\int_{\Omega
}\{\nu k(k-1)\ |\omega |^{k-2}\left\vert \nabla \omega \right\vert
^{2}+\lambda (k-1)|\omega |^{k+1}\}\,dtd\mathbf{x}  \notag \\
& \mathbf{\leqslant }||\omega _{0}||_{L_{k}(\Omega
)}^{k}+\int_{0}^{t_{0}}\int_{\Omega }\lambda (k-1)|\omega |^{k}\ |h|\ dtd%
\mathbf{x+}\int_{0}^{t_{0}}M(\mathbf{v}(t))\int_{\Gamma }\ |b|^{k}\ dtd%
\mathbf{x}.  \label{h10}
\end{align}%
Let $\widetilde{h}:=h-h_{a}$ (see (\ref{haha})). Multiplying the equation of
(\ref{h2}) by $|\widetilde{h}|^{k}sign(\widetilde{h}),$ we obtain 
\begin{equation*}
k\int_{\Omega }|\widetilde{h{}}|^{k-1}|\nabla \widetilde{h{}}|^{2}d\mathbf{x}%
+||\widetilde{h{}}||_{L_{k+1}(\Omega )}^{k+1}\leqslant ||\widetilde{h{}}%
||_{L_{k+1}(\Omega )}^{k}\left\Vert \omega \right\Vert _{L_{k+1}(\Omega )},
\end{equation*}%
that gives%
\begin{equation}
||\widetilde{h{}}||_{L_{k+1}(\Omega )}\leqslant \left\Vert \omega
\right\Vert _{L_{k+1}(\Omega )}.  \label{h9}
\end{equation}%
By the classical estimate for a elliptic problem (\ref{h2}) (see \cite%
{LadyUral68}) 
\begin{equation}
||h||_{W_{p}^{2}(\Omega )}\leqslant Cp(||\omega ||_{L_{p}(\Omega )}+\Vert
a\Vert _{W_{p}^{2}(\Gamma )})\quad \quad \text{for}\quad 2<p<\infty
\label{ell-6}
\end{equation}%
and from the embedding theorem $W_{p}^{1}(\Omega )\hookrightarrow C^{\alpha
}(\bar{\Omega})$ for $\alpha :=1-\frac{2}{p},$ in particular for the
solution $h_{a}$ of the elliptic equation (\ref{haha}) \ we have 
\begin{equation*}
||h_{a}||_{C(\overline{\Omega })}(t)\leqslant C||h_{a}||_{W_{p}^{2}(\Omega
)}(t)\leqslant C\Vert a\Vert _{W_{p}^{2}(\Gamma )}(t)<C.
\end{equation*}%
Therefore Holder's inequality implies%
\begin{equation*}
\int_{\Omega }|h|\ |\omega |^{k}\ d\mathbf{x}\leqslant ||\widetilde{h{}}%
||_{L_{k+1}(\Omega )}\left\Vert \omega \right\Vert _{L_{k+1}(\Omega
)}^{k}+C\left\Vert \omega \right\Vert _{L_{k}(\Omega )}^{k}.
\end{equation*}%
Then from (\ref{h10})-(\ref{h9}), we deduce 
\begin{equation}
\left\Vert \omega \right\Vert _{L_{k}(\Omega )}^{k}(t_{0})\leqslant
C^{k}+\int_{0}^{t_{0}}\left\{ kC\left\Vert \omega \right\Vert _{L_{k}(\Omega
)}^{k}+M(\mathbf{v}(t))\int_{\Gamma }\ |b|^{k}\ d\mathbf{x}\right\} dt.
\label{ha}
\end{equation}%
Using the inequality $(a+b)^{k}\leqslant 2^{k-1}(a^{k}+b^{k})$ with $%
a,b\geqslant 0$ and the regularity assumptions (\ref{reg1})-(\ref{reg2}),
the estimate\ (\ref{ell-6}), the embedding theorem $W_{k}^{1}(\Omega
)\hookrightarrow C(\overline{\Omega }),$ we have that%
\begin{equation*}
M(\mathbf{v}(t))\int_{\Gamma }\ |b|^{k}\ d\mathbf{x}\leqslant
C^{k}(1+\left\Vert \nabla h\right\Vert _{L_{k}(\Gamma )}^{k\varkappa
+1})\leqslant C^{k}(1+\left\Vert \omega \right\Vert _{L_{k}(\Omega )}^{k})
\end{equation*}%
for any$\quad k>k_{0}:=\max (2,1/(1-\varkappa )).$ Hence substituting the
last inequality in the last term of (\ref{ha}), we get 
\begin{equation*}
\left\Vert \omega \right\Vert _{L_{\overline{k}}(\Omega )}^{\overline{k}%
}(t_{0})\leqslant C(1+\int_{0}^{t_{0}}\left\Vert \omega \right\Vert _{L_{%
\overline{k}}(\Omega )}^{\overline{k}}\ dt).
\end{equation*}%
with $\overline{k}:=k_{0}+1.$ Applying Gronwall's inequality we get $%
\left\Vert \omega \right\Vert _{L_{\infty }(0,T;L_{\overline{k}}(\Omega
))}\leqslant C$ \ and $||M(\mathbf{v}(t))||_{L_{\infty }(0,T)},$ $\ \
||b||_{L_{\infty }(0,T;C(\Gamma ))}\leqslant C,$ hence by (\ref{ha}) 
\begin{equation*}
\left\Vert \omega \right\Vert _{L_{k}(\Omega )}^{k}(t_{0})\leqslant
C^{k}+kC\int_{0}^{t_{0}}\left\Vert \omega \right\Vert _{L_{k}(\Omega
)}^{k}dt.
\end{equation*}%
Using Gronwall's inequality and passing to the limit on $k\rightarrow \infty 
$\ we derive (\ref{est-2}).

The estimates (\ref{est3}) follow from (\ref{ell-6}) and the embedding
theorem $W_{p}^{1}(\Omega )\hookrightarrow C^{\alpha }(\bar{\Omega})$ with $%
\alpha :=1-\frac{2}{p}.$

\textbf{2}$^{\mathbf{nd}}$\textbf{\ case:} $g(\omega )=\omega -\omega ^{2}.$
\ Let us deduce the a priori estimate (\ref{wa2}) for the solution $\omega $
of \eqref{h1}-\eqref{h2}. \ Solving \eqref{h1}-\eqref{h2} \ \ we can
consider that $g$ is a \ smooth function being an approximation of $\omega
-\omega ^{2}$ on $[0,1]$ and 
\begin{equation}
g(\omega )=0\quad \quad \text{for}\quad \omega \leqslant 0\quad \text{and}%
\quad \omega \geqslant 1.  \label{const}
\end{equation}%
If we take $\eta (\omega ):=|\omega |_{-}$ \ \ (and $\eta (\omega ):=|\omega
-1|_{+},$ \ respectively) \ in \eqref{ww}, by the assumptions (\ref{wa1}) we
immediately obtain the estimate (\ref{wa2}). As a consequence of it we will
have that the solution of \eqref{sec1eq1}-\eqref{sec1eq6} satisfies (\ref%
{wa2}) \ too.

\bigskip

\subsubsection{\protect\vspace{1pt} Leray-Schauder`s fixed point argument}

Now we are able to construct a compact operator, which fulfills
Leray-Schauder`s fixed point theorem. We choose some function $\widetilde{%
\omega }\in C(\overline{\Omega }_{T}).$ The elliptic problem 
\begin{equation}
\left\{ 
\begin{array}{ll}
-\Delta \widetilde{h}+\widetilde{h}=\widetilde{\omega } & \quad \quad \text{%
in}\quad \Omega _{T}, \\ 
\widetilde{h}=a(t,\mathbf{x}) & \quad \quad \text{on}\quad \Gamma _{T}.%
\end{array}%
\right. \;  \label{o-5}
\end{equation}%
has an unique solution $\widetilde{h}\in L_{\infty }(0,T;C^{1+\alpha }(%
\overline{\Omega }))$ satisfying (\ref{ell-6}). \ Then we consider the
parabolic problem 
\begin{equation}
\left\{ 
\begin{array}{l}
\omega _{t}=\nu \bigtriangleup \omega -\lambda div\left( \mathbf{g}\right)
\quad \quad \text{in}\quad \Omega _{T},\;\;\text{where}\;\mathbf{\ \ }%
\widetilde{\mathbf{v}}:=-\nabla \widetilde{h},\quad \mathbf{g}:=g(\widetilde{%
\omega })\widetilde{\mathbf{v}}, \\ 
\\ 
\nu \frac{\partial \omega }{\partial \mathbf{n}}+\lambda M(\mathbf{v}%
(t))(\omega -b(\cdot ,\text{\textbf{$\cdot $}},\frac{\partial \widetilde{h}}{%
\partial \mathbf{n}}))=0\quad \quad \text{\ on}\quad \text{\ }\Gamma _{T},
\\ 
\\ 
\omega |_{t=0}=\lambda \omega _{0}\quad \quad \text{\ in}\quad \text{\ }%
\Omega .%
\end{array}%
\right.  \label{o-7}
\end{equation}%
We have that $\left\Vert \mathbf{g}\right\Vert _{L^{\infty }(\Omega
_{T})}\leqslant C\left\Vert \widetilde{\omega }\right\Vert _{C(\overline{%
\Omega }_{T})}^{2}.$ Applying Theorem 5.1, p.170, Theorem 6.1, p. 178 and
Theorem 10.1, p. 204 of \cite{LadySolonUral68}, we get that the system (\ref%
{o-7}) has an unique solution $\omega \in C^{\alpha /2,\alpha }(\overline{%
\Omega }_{T})\ \;\;\;$for some $\;\alpha \in (0,1),$ such that%
\begin{equation}
\left\Vert \omega \right\Vert _{C^{\alpha /2,\alpha }(\overline{\Omega }%
_{T})}\leqslant \;C\mathbf{(}\nu \mathbf{)}\;(||b(\cdot ,\cdot ,\frac{%
\partial \widetilde{h}}{\partial \mathbf{n}}(\cdot ,\cdot ))||_{C^{1,1}(%
\overline{\Omega }_{T})}+\left\Vert \mathbf{g}\right\Vert _{L^{\infty
}(\Omega _{T})}),  \label{ch1}
\end{equation}%
with the constant $\ C\mathbf{(}\nu \mathbf{)}$ depending on $\ \nu .$
Therefore we construct the operator 
\begin{equation}
\omega :=B\left[ \widetilde{\omega },\lambda \right] ,  \label{o-8}
\end{equation}%
such that $B\left[ \widetilde{\omega },0\right] =0.$ By (\ref{ch1}) and
Theorem 4.5, p.166 of \ \cite{LadySolonUral68} $\ B\ $is a \textit{compact} 
\textit{continuous} operator of \ the Banach space $C(\overline{\Omega }%
_{T})\times \lbrack 0,1]\ $into $C(\overline{\Omega }_{T}).$ \ Let us note
that the estimate (\ref{est-2})\ is valid for all $\omega =B\left[ \omega
,\lambda \right] ,$ being a solution of (\ref{h1})-(\ref{h2}). Therefore $B\ 
$fulfills all conditions of \ Leray-Schauder's \ fixed point theorem and as
a consequence of it\textbf{\ }the problem \eqref{pe1}-\eqref{pe2} is
solvable in the Banach space $C(\overline{\Omega }_{T}).$ Taking into
account (\ref{ch1}) and applying classical results of \cite{LadyUral68} for
elliptic problem (\ref{pe2}),$\ $we obtain that\ $h\in C^{\alpha /2,2+\alpha
}(\overline{\Omega }_{T}).$\ Finally with the help of Theorem 12.2, p. 224, 
\cite{LadySolonUral68}, we derive that $\omega \in C^{1+\alpha /2,2+\alpha }(%
\overline{\Omega }_{T}).$

\subsubsection{\protect\vspace{1pt} Proof of a priori estimates (\protect\ref%
{est22})-(\protect\ref{est2})}

The estimate (\ref{est22}) follows from (\ref{est-2})-(\ref{est3}) and (\ref%
{h10}) for $k=2$.

In view of (\ref{pe2}) the function $H:=\frac{\partial (h-h_{a})}{\partial t}%
\ $solves the following elliptic problem 
\begin{equation*}
\left\{ 
\begin{array}{l}
-\Delta H+H=\mathrm{div}\,\mathbf{G}\;\text{ in }\;\Omega , \\ 
H\big|_{\Gamma }=0,%
\end{array}%
\right.
\end{equation*}%
with $\mathbf{G}:=\nu \nabla \omega -\mathbf{v}g(\omega ),$ such that $||%
\mathbf{G}||_{L_{2}(\Omega _{T})}\leqslant C.$ Applying the classical
results of \cite{LadyUral68}, we derive 
\begin{equation*}
\left\Vert \nabla H(t)\right\Vert _{L_{2}(\Omega )}\leqslant C\left\Vert 
\mathbf{G}(t)\right\Vert _{L_{2}(\Omega )}\;\text{ for }\;t\in (0,T),
\end{equation*}%
\bigskip which gives the estimate (\ref{est1}). \ \ 

By virtue of \eqref{est3} and \eqref{est1}, we have \ that 
\begin{equation*}
\nabla (h-h_{a})\in \mathcal{\Re }\left[ \Omega \right] :=L_{2}\big(%
0,T;W_{2}^{1}(\Omega )\big)\cap W_{2}^{1}\big(0,T;L_{2}(\Omega )\big),
\end{equation*}%
such that 
\begin{equation*}
||\nabla (h-h_{a})||_{\mathcal{\Re }\left[ \Omega \right] }\leqslant C.
\end{equation*}%
Let us choose a neighborhood $G$ of $\Gamma $ and an orthogonal coordinate
system $(y_{1},y_{2}),$ such that $G=\{\mathbf{y}=(y_{1},y_{2}):(y_{1},0)\in
\Gamma ,$ $y_{2}\in \lbrack 0,\eta ]\}$ for some $\eta >0$ and this new
coordinate system coincides with the tangential-normal coordinate system
along $\Gamma .$ The space $\mathcal{\Re }\left[ G\right] $ can be
represented as 
\begin{align*}
\mathcal{\Re }\left[ G\right] =\Big\{\,\phi (y_{2},t,y_{1})\in L_{2}(0,\eta
;\;L_{2}(0,T;& W_{2}^{1}(\Gamma )))\cap W_{2}^{1}(0,T;L_{2}(\Gamma )): \\
& \partial _{y_{2}}\,\phi \in L_{2}\big(0,\eta ;L_{2}(0,T;L_{2}(\Gamma ))%
\big)\Big\}.
\end{align*}%
By virtue of the theorem 3.2 of \cite{LionsMag}, the mapping $A:\,\phi
\rightarrow \phi |_{y_{2}=0}$ is well-defined on $\mathcal{\Re }\left[ G%
\right] .$ Furthermore, the operator $\,\,\phi \rightarrow A(\,\phi )$ from $%
\mathcal{\Re }\left[ G\right] $ to 
\begin{equation*}
\Big[L_{2}\big(0,T;W_{2}^{1}(\Gamma )\big)\cap W_{2}^{1}\big(%
0,T;L_{2}(\Gamma )\big),L_{2}\big(0,T;L_{2}(\Gamma )\big)\Big]_{\frac{1}{2}%
}\equiv \mathcal{P}\left[ \Gamma \right]
\end{equation*}%
is continuous and surjective. Here $[X,Y]_{\varepsilon },$ $\varepsilon \in
\lbrack 0,1],$ are the intermediate spaces of \ Banach spaces $X$ and $Y$ as
defined in \cite{LionsMag}. Therefore the value of $\frac{\partial (h-h_{a})%
}{\partial \mathbf{n}}|_{\Gamma }$ belongs to the space $\mathcal{P}\left[
\Gamma \right] $ and satisfies the estimate (\ref{est2}).

\subsection{\protect\vspace{1pt}Proof of Lemma \ \protect\ref{lemma2}}

\label{sec 4.2}

Let us define the vector function $\mathbf{F}_{\alpha }=\mathbf{F}_{\alpha
}(t,\mathbf{x})$, having the following components%
\begin{equation*}
F_{1,\alpha }:=\int\limits_{\mathbb{R}}f(\xi ,\cdot ,\cdot )\,\alpha (\xi
)\,d\xi \;\text{ and }\;F_{2,\alpha }:=\mathbf{v}\int\limits_{\mathbb{R}%
}\,g^{\prime }(\xi )\,f(\xi ,\cdot ,\cdot )\,\alpha (\xi )\,d\xi
\end{equation*}%
for a fixed $\alpha \in C_{0}^{\infty }(\mathbb{R}).$ From \eqref{eq2.9} we
have that in the distributional sense 
\begin{eqnarray*}
\mathrm{div}_{t,\mathbf{x}}\mathbf{F}_{\alpha } &=&\int\limits_{\mathbb{R}%
}(\partial _{t}f+\mathrm{div}_{\mathbf{x}}(g^{\prime }(\xi )\mathbf{v}%
f))\alpha (\xi )\,d\xi \\
&=&\int\limits_{\mathbb{R}}\left\{ -m(\xi ,\cdot ,\cdot )+g(\xi )(fh-\rho
)\right\} \alpha ^{\prime }(\xi )\,d\xi \in \mathcal{M}(\Omega _{T}).
\end{eqnarray*}%
Therefore by Theorem 2.1 of \cite{chen} there exists a continuous linear
functional $\mathbf{F}_{\alpha }\cdot \mathbf{n}_{(t,\mathbf{x})}$ \ over $%
L_{\infty }(\Sigma ).$ Here $\Sigma $ is the boundary of the domain $\ \ 
\overline{\Omega }_{T}$ \ and $\ \mathbf{n}_{(t,\mathbf{x})}$ is the
external normal to $\Sigma .$

\subsubsection{Proof of the part 1) of Lemma \ \protect\ref{lemma2}}

In particular, if we take an arbitrary $\psi \in C^{\infty }([0,T]\times
\Omega )$ with a compact support on $\Omega ,$ then 
\begin{eqnarray*}
<\mathbf{F}_{\alpha }\cdot \mathbf{n}_{(t,\mathbf{x})}\Big|_{t=0},\psi
>&=&\lim_{\varepsilon \rightarrow 0}\frac{1}{\varepsilon }%
\int\limits_{0}^{\varepsilon }\int\limits_{\Omega }[\int\limits_{\mathbb{R}%
}f(\xi ,t,\mathbf{x})\alpha (\xi )\,d\xi ]\psi (t,\mathbf{x})\,\,dt\,d%
\mathbf{x} \\
&=&\lim_{\varepsilon \rightarrow 0}\int\limits_{\mathbb{R}\times \Omega }[%
\frac{1}{\varepsilon }\int\limits_{0}^{\varepsilon }f(\xi ,t,\mathbf{x}%
)\,dt]\alpha (\xi )\psi (0,\mathbf{x})\,\,d\xi \,d\mathbf{x}.
\end{eqnarray*}%
Using the dominated convergence theorem and \ the sequence $\left\{ \frac{1}{%
\varepsilon }\int_{0}^{\varepsilon }f(\cdot ,t,\cdot )\,dt\right\}
_{\varepsilon \in (0,T)}$ is bounded by $0$ from below and $1$ from above,
we get that $\mathbf{F}_{\alpha }\cdot \mathbf{n}_{(t,\mathbf{x})}|_{t=0}$
is a bounded function in $L^{\infty }(\mathbb{R}\times \Omega ),$ which is
equal to $f^{0}\alpha .$ The function $f^{0}$ satisfies \eqref{2.10}, as a
consequence of $0\leqslant f\leqslant 1,$ $\ f_{\xi }^{\prime }\leqslant 0$
(see \eqref{eq2.88}).

Next, let us choose $\varphi :=O_{\varepsilon }(t)\phi \,$ in \eqref{2.60}
with a non-negative$\ \phi \in C_{0}^{\infty }(\Omega )$ and 
\begin{equation*}
O_{\varepsilon }(t):=\left\{ 
\begin{array}{c}
\frac{\varepsilon -t}{\varepsilon }\text{ \ \ if }0\leqslant t\leqslant
\varepsilon , \\ 
0\text{ \ \ if }t\in (\varepsilon ,T].%
\end{array}%
\right.
\end{equation*}
\ Taking $\varepsilon \rightarrow 0$ in an obtained inequality, we derive 
\begin{equation*}
-\int\limits_{\Omega }\int\limits_{\xi }^{+\infty }f^{0}(s,\mathbf{x})\,\phi
(\mathbf{x})\,ds\,d\mathbf{x}+\int\limits_{\Omega }|\omega _{0}-\xi
|_{+}\,\phi (\mathbf{x})\,d\mathbf{x}\geqslant 0\text{ \ a.e. in }\mathbb{R}%
\times \Omega .
\end{equation*}%
Therefore the function $m_{+}^{0}(\xi ,\mathbf{x}):=|\omega _{0}(\mathbf{x}%
)-\xi |_{+}-\int\limits_{\xi }^{+\infty }f^{0}(s,\mathbf{x})\,ds$ \ in $%
\mathbb{R}\times \Omega $ \ satisfies the first line of \eqref{2.13}.

By the same way if we put the above chosen $\varphi $ \ in \eqref{2.61} and
pass on $\varepsilon \rightarrow 0$ in an obtained inequality, we get 
\begin{equation*}
-\int\limits_{-\infty }^{\xi }(1-f^{0}(s,\mathbf{x}))\,ds+|\omega _{0}-\xi
|_{-}d\mathbf{x}\geqslant 0\text{ \ a.e. in }\mathbb{R}\times \Omega .
\end{equation*}
Hence the function $m_{-}^{0}(\xi ,\mathbf{x}):=|\omega _{0}-\xi
|_{-}-\int\limits_{-\infty }^{\xi }(1-f^{0}(s,\mathbf{x}))\,ds$ in $\mathbb{R%
}\times \Omega $ satisfies the second line of \eqref{2.13}.

\subsubsection{Proof of the part 2) of Lemma \ \protect\ref{lemma2}}

\bigskip \label{sec 333}

Now we take an arbitrary $\psi \in C^{\infty }((0,T)\times \overline{\Omega }%
)$ with a compact support on $(0,T)$\ and using $\mathbf{v}\in L_{\infty
}(0,T,C^{\alpha }(\overline{\Omega })),$\ we have 
\begin{eqnarray*}
<\mathbf{F}_{\alpha }\cdot \mathbf{n}_{(t,\mathbf{x})}\Big|_{\Gamma
_{T}},\psi > &=&\lim_{\varepsilon \rightarrow 0}\frac{1}{\varepsilon }%
\int\limits_{0}^{\varepsilon }\int\limits_{\Gamma _{T}}[\int\limits_{\mathbb{%
R}}g^{\prime }(\xi )(\mathbf{v}(t,\mathbf{x}^{\prime }\mathbf{)}\cdot 
\mathbf{n}(\mathbf{x)})f(\xi ,t,\mathbf{x}^{\prime })\alpha (\xi )\,d\xi
]\psi (t,\mathbf{x}^{\prime })\,dtd\mathbf{x}\,ds\qquad \\
&=&\int\limits_{\mathbb{R}\times \Gamma _{T}}g^{\prime }(\xi )\,\alpha (\xi
)\,(\mathbf{v}\cdot \mathbf{n})|_{\Gamma }\,\psi |_{\Gamma
}\lim_{\varepsilon \rightarrow 0}[\frac{1}{\varepsilon }\int\limits_{0}^{%
\varepsilon }f(\xi ,t,\mathbf{x}^{\prime })\,ds]\,d\xi dtd\mathbf{x,}
\end{eqnarray*}%
here $\mathbf{x}^{\prime }:\mathbf{=x}-s\,\mathbf{n}(\mathbf{x}).$ Since the
sequence $\left\{ \frac{1}{\varepsilon }\int_{0}^{\varepsilon }f(\cdot
,\cdot ,\mathbf{x}^{\prime })\,ds\right\} _{\varepsilon \in (0,\varepsilon
_{0})}$ is bounded by $0$ from below and $1$ from above, \ then $\mathbf{F}%
_{\alpha }\cdot \mathbf{n}_{(t,\mathbf{x)}}|_{\Gamma _{T}}$ is a bounded
function in $L^{\infty }(\mathbb{R}\times \Gamma _{T}),$ which is equal to $%
f^{\Gamma }\alpha .$ The function $f^{\Gamma }$ satisfies \eqref{eq2.15}, as
a consequence of $0\leqslant f\leqslant 1,$ $\ $ $f_{\xi }^{\prime
}\leqslant 0.$

Let 
\begin{equation}
{\mathbf{0}}_{\varepsilon }(\mathbf{x}):=\left\{ 
\begin{array}{l}
\frac{\varepsilon -d}{\varepsilon },\quad \mbox{ if }0\leqslant d(\mathbf{x}%
)<\varepsilon , \\ 
0,\quad \quad \,\mbox{ if
}\mathbf{x}\in \Omega \mbox{ and }d(\mathbf{x})\geqslant \varepsilon .%
\end{array}%
\right.  \label{xunit}
\end{equation}%
be the approximation of the zero function on $\overline{\Omega }.$ Here $d(%
\mathbf{x}):=\min_{\mathbf{y}\in \Gamma }|\mathbf{x}-\mathbf{y}|$ is the
distance between $\mathbf{x}\in \overline{\Omega }$ and the boundary $\Gamma
.$

Let us choose in \eqref{2.60} $\varphi :={\mathbf{0}}_{\varepsilon }(\mathbf{%
x})\,\phi $ for a non-negative $\phi \in C^{\infty }((0,T)\times \overline{%
\Omega }),$ having a compact support on $(0,T).$ Taking $\varepsilon
\rightarrow 0$ in a deduced inequality, since $\ \nabla d|_{\Gamma }=-%
\mathbf{n}$ we derive 
\begin{equation*}
\int\limits_{\Gamma _{T}}\int\limits_{\xi }^{+\infty }g^{\prime }(s)\,(%
\mathbf{v}\cdot \mathbf{n})\,\,f^{\Gamma }(s,t,\mathbf{x})\,ds\,\,\phi (t,%
\mathbf{x})\,dtd\mathbf{x}+M\int\limits_{\Gamma _{T}}|b-\xi |_{+}\,\phi (t,%
\mathbf{x})\,dtd\mathbf{x}\geqslant 0.
\end{equation*}%
Therefore the function $m_{+}^{\Gamma }(\xi ,t,\mathbf{x}):=M\,\,|b(t,%
\mathbf{x},\frac{\partial h}{\partial \mathbf{n}})-\xi
|_{+}+\int\limits_{\xi }^{+\infty }g^{\prime }(s)\,(\mathbf{v}\cdot \mathbf{n%
})\,\,f^{\Gamma }(s,t,\mathbf{x})\,ds$ on \ $\mathbb{R}\times \Gamma _{T}$
satisfies the first two lines of \eqref{2.16}.

By the same way if we put the above chosen $\varphi $ \ in \eqref{2.61} and
pass on $\varepsilon \rightarrow 0$ in a deduced inequality, we derive 
\begin{equation*}
\int\limits_{-\infty }^{\xi }g^{\prime }(\xi )(\mathbf{v}\cdot \mathbf{n}%
)(1-f^{\Gamma }(s,t,\mathbf{x}))\,ds+M|b-\xi |_{-}\geqslant 0,
\end{equation*}%
Hence the function $m_{-}^{\Gamma }(\xi ,t,\mathbf{x}):=M|b(t,\mathbf{x},%
\frac{\partial h}{\partial \mathbf{n}})-\xi |_{-}+\int\limits_{-\infty
}^{\xi }g^{\prime }(s)(\mathbf{v}\cdot \mathbf{n})(1-f^{\Gamma }(s,t,\mathbf{%
x}))\,ds$ on \ $\mathbb{R}\times \Gamma _{T}$ satisfies the last two lines
of \eqref{2.16}.

\bigskip

\subsection{Proof of Lemma \protect\ref{lemma33}}

\label{sec 33}

The existence of the traces for $F$ \ both on the boundary $\Gamma $ and for
a.e. $t\in \lbrack 0,T]$\ do not follow from the \textit{inequality} (\ref%
{eqq1}) (in contrast with $f$ which satisfies the equality (\ref{2.8})),
that presents the main difficulty in the deduction of (\ref{uha}). To show
it a simple inequality 
\begin{equation}
-\frac{1}{\varepsilon }\int_{0}^{\varepsilon }z^{2}(s)ds\leqslant -\left( 
\frac{1}{\varepsilon }\int_{0}^{\varepsilon }z(s)ds\right) ^{2}
\label{holder}
\end{equation}%
will be used, which is valid for any positive integrable function $z=z(s).$

Let us introduce the functions 
\begin{equation}
{\mathbf{1}_{\varepsilon }}(\mathbf{x}):=\left\{ 
\begin{array}{l}
\frac{d(\mathbf{x})}{\varepsilon },\quad \mbox{ if }0\leqslant d(\mathbf{x}%
)\leqslant \varepsilon , \\ 
1,\quad \quad \mbox{ if
}\mathbf{x}\in \Omega \mbox{ and }\varepsilon <d(\mathbf{x}),%
\end{array}%
\right. \quad {{1}_{\varepsilon }}(t):=\left\{ 
\begin{array}{l}
0,\quad \mbox{ if }t<0, \\ 
\frac{t}{\varepsilon },\quad \mbox{ if }0\leqslant t\leqslant \varepsilon ,
\\ 
1,\quad \mbox{ if }\sigma <t.%
\end{array}%
\right.  \label{222}
\end{equation}%
The distance function $d=d(\mathbf{x})$ is defined in the section \ref{sec
333}. If we choose 
\begin{equation*}
\varphi (\xi ,t,\mathbf{x}):=(1{_{\varepsilon _{0}}}(\xi +\varepsilon
_{0}^{-1})-1{_{\varepsilon _{0}}}(\xi -\varepsilon _{0}^{-1}))(1{%
_{\varepsilon }}(t)-1{_{\varepsilon }}(t-t_{0}+\varepsilon ))\mathbf{1}%
_{\varepsilon }(\mathbf{x})
\end{equation*}%
with $t_{0}\in (2\varepsilon ,T)$ as the test function \ in the respective
integral formulation of \eqref{eqq1} \ and pass to the limit on $\varepsilon
_{0}\rightarrow 0,$ we derive 
\begin{align}
\frac{1}{\varepsilon }\int_{t_{0}-\varepsilon }^{t_{0}}& \left( \int\limits_{%
\mathbb{R}\times \Omega }F\ \mathbf{1}_{\varepsilon }(\mathbf{x})\,d\xi d%
\mathbf{x}\right) dt\leqslant \frac{1}{\varepsilon }\int_{0}^{\varepsilon
}\left( \int_{{\mathbb{R}}\times \Omega }F\ \mathbf{1}_{\varepsilon }(%
\mathbf{x})\ d\xi d\mathbf{x}\right) \,dt  \notag \\
& +\int_{{\mathbb{R\times }}[0,T]}\left( \frac{1}{\varepsilon }%
\int_{0\leqslant d(\mathbf{x})\leqslant \varepsilon }|g^{\prime }(\xi )(%
\mathbf{v\cdot }\nabla d)|_{-}\ F\ d\mathbf{x}\right) \,d\xi dt  \notag \\
& +\int\limits_{0}^{t_{0}}\left( \int\limits_{\mathbb{R}\times \Omega
}\left\vert \left( g(\xi )(\omega -\xi )\right) _{\xi }^{\prime }F+g^{\prime
}(\xi )(1-2f)[\rho -f\omega ]\right\vert \ \mathbf{1}_{\varepsilon }(\mathbf{%
x})\ d\xi d\mathbf{x}\right) \,dt  \notag \\
& =:A_{1}^{\varepsilon }+A_{2}^{\varepsilon }+A_{3}^{\varepsilon },
\label{ddF}
\end{align}%
where we have used that $F=0$ and $\rho -f\omega =0$ for $|\xi |>R.$

From (\ref{holder}) it follows that 
\begin{equation*}
A_{1}^{\varepsilon }\leqslant \int_{{\mathbb{R}}\times \Omega }\left[ \left( 
\frac{1}{\varepsilon }\int_{0}^{\varepsilon }f(t)dt\right) -\left( \frac{1}{%
\varepsilon }\int_{0}^{\varepsilon }f(t)dt\right) ^{2}\right] \ \mathbf{1}%
_{\varepsilon }(\mathbf{x})\ d\xi d\mathbf{x}.
\end{equation*}%
The sequence $\frac{1}{\varepsilon }\int_{0}^{\varepsilon }f(t)dt$ is
uniformly bounded on $\varepsilon $\ by $0$ and $1,$ converging to $f^{0}$
a.e. on ${\mathbb{R}}\times \Omega $ by Lemma \ref{lemma2}, then in view of
the dominated convergence theorem, we derive%
\begin{equation}
\limsup_{\varepsilon \rightarrow 0}A_{1}^{\varepsilon }\leqslant \int_{{%
\mathbb{R}}\times \Omega }(f^{0}-f{^{0}}^{2})\ d\xi d\mathbf{x}.  \label{oi1}
\end{equation}%
Now we consider the term $A_{2}^{\varepsilon }.$ For a small enough $%
\varepsilon $ there exists an unique projection $\mathbf{x}_{0}=\mathbf{x}%
_{0}(\mathbf{x})$ on the boundary $\Gamma $\ of any point $\mathbf{x}\in
\Omega :\ d(\mathbf{x})<\varepsilon .$ Let us parametrize all points $%
\mathbf{x}_{0}\in \Gamma $ by the length parameter $l\in \lbrack 0,L]$ \
with $L$ being the length of $\ \Gamma .$ \ Let us make the change of
variables $\mathbf{x}=(x_{1},x_{2})\leftrightarrow (l,s)$ in the term $%
A_{2}^{\varepsilon },$ where $s$ is define as $s:=d(\mathbf{x).}$ We have
the following\ \ properties:

1) $\nabla d(\mathbf{x})=-\mathbf{n}(\mathbf{x}_{0})+O(\varepsilon )$ and
the Jacobian $\frac{D(x_{1},x_{2})}{D(l,s)}=1+O(\varepsilon ),$ since $%
\Gamma \in C^{2}$ and $(l,s)$ forms the orthogonal coordinate system at $%
s=0; $

2) $\mathbf{v}(t,\mathbf{x})=$ $\mathbf{v}(t,\mathbf{x}_{0})+O(\varepsilon
^{\alpha }),$ since $\mathbf{v}(t,\cdot )\in C^{\alpha }(\overline{\Omega })$
with $\alpha =1-\frac{2}{p}.$

Therefore using (\ref{holder}) it follows 
\begin{eqnarray*}
A_{2}^{\varepsilon } &\leqslant &\frac{1}{\varepsilon }\int_{0}^{\varepsilon
}\int_{0}^{L}\left[ \int_{{\mathbb{R}}\times \lbrack 0,T]}|g^{\prime }(\xi )(%
\mathbf{v}\cdot \mathbf{n})(t,\mathbf{x}_{0})|_{-}\ F(\xi ,t,\mathbf{x}_{0}-s%
\mathbf{n}(\mathbf{x}_{0}))\ dtd\xi \right] \ dsdl\ +O(\varepsilon ^{\alpha
}) \\
&\leqslant &\int_{{\mathbb{R}}\times \Gamma _{T}}|g^{\prime }(\xi )(\mathbf{v%
}\cdot \mathbf{n})(t,\mathbf{x}_{0})|_{-}\ [f^{\varepsilon }-\left(
f^{\varepsilon }\right) ^{2}]\ d\xi dtd\mathbf{x}_{0}+O(\varepsilon ^{\alpha
})
\end{eqnarray*}%
with $\mathbf{x}_{0}=\mathbf{x}_{0}(l)$ in the 1$^{st}$ inequality and in
the 2$^{nd}$ inequality we denote by $f^{\varepsilon }:=\frac{1}{\varepsilon 
}\int_{0}^{\varepsilon }f(\cdot ,\cdot ,\mathbf{x}_{0}-s\mathbf{n}(\mathbf{x}%
_{0}))ds.$\ Since the sequence $\frac{1}{\varepsilon }\int_{0}^{\varepsilon
}f(\cdot ,\cdot ,\mathbf{x}_{0}-s\mathbf{n}(\mathbf{x}_{0}))ds$ is uniformly
bounded on $\varepsilon $\ by $1$ and converges to $f^{\Gamma }$ a.e. on ${%
\mathbb{R}}\times \Gamma _{T},\ \ $where$\ \ g^{\prime }(\xi )(\mathbf{v}%
\cdot \mathbf{n})\neq 0$ by Lemma \ref{lemma2}, we obtain%
\begin{equation}
\limsup_{\varepsilon \rightarrow 0}A_{2}^{\varepsilon }\leqslant \int_{{%
\mathbb{R}}\times \Gamma _{T}}|g^{\prime }(\xi )(\mathbf{v}\cdot \mathbf{n}%
)(t,\mathbf{x})|_{-}(f^{\Gamma }-\left( {f}^{\Gamma }\right) ^{2})\ d\xi dtd%
\mathbf{x}.  \label{oi2}
\end{equation}%
Using the boundedness of $\omega ,\ g(\xi ),\ g^{\prime }(\xi ),$ $(1-2f)$
in $L_{\infty }$-spaces and \eqref{eq2.21}, we get that the term $%
A_{3}^{\varepsilon }$ is estimated as%
\begin{equation}
A_{3}^{\varepsilon }\leqslant C\int\limits_{0}^{t_{0}}\left( \int\limits_{%
\mathbb{R}\times \Omega }F\ \mathbf{1}_{\varepsilon }(\mathbf{x})\ \,d\xi d%
\mathbf{x}\right) \ dt\mathbf{,}  \label{oi3}
\end{equation}%
with the constant $C$ being \ independent\textit{\ }of $\varepsilon .$
Finally substituting the inequalities (\ref{oi1})-(\ref{oi3}) into (\ref{ddF}%
), taking the limit on $\varepsilon \rightarrow 0$ in the derived inequality
and applying Lebesgue's differentiation theorem to the left part of this
inequality, we get 
\begin{eqnarray*}
\int_{\mathbb{R}\times \Omega }F\ d\xi d{\mathbf{x}}(t_{0}) &\leqslant
&\int_{\mathbb{R}\times \Omega }(f^{0}-\left( f^{0}\right) ^{2})\ d\xi d{%
\mathbf{x}} \\
&+&\int\limits_{\mathbb{R}\times \Gamma _{T}}{|g^{\prime }(\xi )(\mathbf{v}%
\cdot \mathbf{n})|}_{-}\,\,(f^{\Gamma }-\left( f^{\Gamma }\right)
^{2})\,\,d\xi dtd\mathbf{x}+C\int\limits_{0}^{t_{0}}\left( \int\limits_{%
\mathbb{R}\times \Omega }F\ \,d\xi d\mathbf{x}\right) \ dt,
\end{eqnarray*}%
that implies (\ref{uha}) by Gronwall's inequality. $\blacksquare $

\subsection{Proof of Lemma \protect\ref{lemma3}}

\label{sec 4.4} \setcounter{equation}{0}

By virtue of Lemma \ref{lemma2}, 1) and using integration on $\xi $ by parts
we have that 
\begin{eqnarray*}
0 &\geqslant &\int\limits_{\mathbb{R}}f^{0}(f^{0}-1)\,d\xi
=\int\limits_{-\infty }^{\omega _{0}}f^{0}\left\{ \mathrm{sign}_{+}(\omega
_{0}-\xi )-1+\partial _{\xi }m_{-}^{0}\right\} \,d\xi \\
&+&\int\limits_{\omega _{0}}^{+\infty }\left\{ \mathrm{sign}_{+}(\omega
_{0}-\xi )+\partial _{\xi }m_{+}^{0}\right\} (f^{0}-1)\,d\xi
=\int\limits_{-\infty }^{\omega _{0}}f^{0}\,\partial _{\xi }m_{-}^{0}d\xi
+\int\limits_{\omega _{0}}^{+\infty }\partial _{\xi
}m_{+}^{0}\,(f^{0}-1)\,d\xi \\
&=&(f^{0}m_{-}^{0})|_{\xi =\omega _{0}}-\int\limits_{-\infty }^{\omega
_{0}}\partial _{\xi }f^{0}\,m_{-}^{0}\,d\xi +m_{+}^{0}(1-f^{0})|_{\xi
=\omega _{0}}-\int\limits_{\omega _{0}}^{+\infty }m_{+}^{0}\,\partial _{\xi
}f^{0}\,d\xi \geqslant 0
\end{eqnarray*}%
a.e. \ on \ $\Omega .$ Let us note that the justification of the last
identity (which is formal) can be done by a standard procedure: mollifying \
the function $f^{0}$ and passing to the limit (see, for instance, \ \cite%
{lellis}).

Taking into account that $s=-|s|_{-}$ for $s<0,$\ then by Lemma \ref{lemma2}%
, 2) and the integration on $\xi $ by parts, we have%
\begin{align*}
0& \leqslant -\int\limits_{\mathbb{R}}|g^{\prime }(\xi )(\mathbf{v}\cdot 
\mathbf{n})|_{-}\,f^{\Gamma }(f^{\Gamma }-1)\,d\xi =\int\limits_{-\infty
}^{b}f^{\Gamma }\{M\left[ 1-\mathrm{sign}_{+}(b-\xi )\right] -\partial _{\xi
}m_{-}^{\Gamma })\}\,d\xi \\
& +\int\limits_{b}^{+\infty }[-M\mathrm{sign}_{+}(b-\xi )-\partial _{\xi
}m_{+}^{\Gamma }](f^{\Gamma }-1)\,d\xi =-\left[ (f^{\Gamma }m_{-}^{\Gamma
})|_{\xi =b-0}-\int\limits_{-\infty }^{b}\partial _{\xi }f^{\Gamma
}m_{-}^{\Gamma }\,d\xi \right. \\
& +\left. m_{+}^{\Gamma }(1-f^{\Gamma })|_{\xi
=b+0}-\int\limits_{b}^{+\infty }m_{+}^{\Gamma }\partial _{\xi }f^{\Gamma
}d\xi \right] \leqslant 0\text{\qquad a.e. on }\Gamma _{T}.
\end{align*}

\section{Acknowledgement}

N.V. Chemetov thanks the support from FCT and Project "F\'{\i}sica - Matem%
\'{a}tica", PTDC / MAT / 69635/ 2006, financed by FCT. \ 

N.V. Chemetov thanks the support from Program "Conv\'{e}nio GRICES / CAPES",
financed by FCT, project \textquotedblleft Euler Equations and related
problems\textquotedblright , cooperation between Portugal (Universidade de
Lisboa) and Brasil (Universidade Estadual de Campinas).

\end{document}